%% file: a5.tex
\documentclass{article}
\usepackage{amsmath,amsfonts,epsf,chicago,pb-diagram,lamsarrow,pb-lams}

\newcommand{\rf}[1]{(\expandafter\ref{#1})}
\newcommand{\ct}[1]{\citeANP{#1}~[\citeyearNP{#1}]}
\newcommand{\ctt}[2]{\citeANP{#1}~[\citeyearNP{#1}, \citeyearNP{#2}]}
\newcommand{\lb}[1]{\expandafter\label{#1}}

\newcommand{\onm}\operatorname
\newcommand{\Id}{{\mathchoice
 {\mbox{Id}}{\mbox{Id}}{\mbox{\scriptsize Id}}{\mbox{\scriptsize Id}}}}
\newcommand{\bset}[2]{\bigl\{\,#1\bigm|#2\,\bigr\}}

\newcommand{\dbyd}[2]{{\displaystyle\left.\frac d{d#1}\right|_{#1=#2}}}

\newcommand{\Wred}{W_{\!\mathrm{red}}}

\newcommand{\Hred}{H_{\mathrm{red}}}
\newcommand{\Hdrift}{H_{\mathrm{drift}}}
\newcommand{\Hnil}{H_{\mathrm{nil}}}
\newcommand{\Pdrift}{P_{\mathrm{drift}}}
\newcommand{\Hrsd}{H_{\mathrm{rsd}}}
\newcommand{\Hrmd}{H_{\mathrm{rmd}}}
\newcommand{\Hxie}{{H_{\xi_e}}}
\newcommand{\xie}{{\xi_e}}
\newcommand{\mue}{{\mu_e}}

\newcommand{\defemph}[1]{{\slshape\sffamily#1}}

\begin{document}\allowdisplaybreaks

\title{
Dynamics near Relative Equilibria: Nongeneric Momenta at a 1:1
Group--Reduced Resonance}
\author{G. W. Patrick\\
Department of Mathematics and Statistics\\
University of Saskatchewan\\Saskatoon, Saskatchewan, 
S7N~5E6\\Canada}
\date{February, 1999\\
{\small\em To appear in  Mathematische Zeitschrift}}
\maketitle

\begin{abstract}
An interesting situation occurs  when the  linearized dynamics of  the
shape   of  a  formally stable  Hamiltonian    relative equilibrium at
nongeneric momentum  1:1~resonates  with a  frequency  of the relative
equilibrium's generator.  In this   case some of the  shape  variables
couple   to  the  group   variables to   first  order  in the momentum
perturbation, and the first order perturbation theory implies that the
relative equilibrium slowly changes orientation in the same way that a
charged particle~\emph{with magnetic moment} moves  on a sphere  under
the influence of a radial magnetic monopole. In  the course of showing
this  a normal form    is constructed for  linearizations of  relative
equilibria  and  for  Hamiltonians   near  group  orbits of   relative
equilibria.
\end{abstract}

{\footnotesize\tableofcontents}

\section{Introduction}

Suppose we are given a  symplectic manifold $(P,\omega)$ and a compact
Lie   group $G$  acting  freely     and symplectically  on  $P$   with
$\mbox{Ad}^*$-equivariant  momentum  mapping  $J:P\rightarrow\mathfrak
g^*$.   Let   the  Lie   algebra  of   $G$   be  $\mathfrak    g$ and
$H:P\rightarrow\mathbb R$ be a $G$-invariant  Hamiltonian.  Let us fix
our    attention on   a   specific  relative  equilibrium~$p_e$,  with
generator~$\xi_e$, so that the integral curve of  the Hamiltonian
vector field $X_H$ at $p_e$ is  $\exp(\xie t)p_e$.  Let~$J(p_e)=\mue$,
let  $G_\mue$  be the isotropy  group  of $\mue$   under the coadjoint
action of  $G$ on $\mathfrak g^*$,  and let  $\mathfrak g_\mue$ be the
Lie  algebra  of $G_\mue$.    For   more information  on   these basic
definitions,                          see,                         for
example,~\ct{AbrahamRMarsdenJE-1978.1},~\ct{MarsdenJE-1992.1},
or~\ct{MarsdenJE-1994.1}.

Suppose that $G=SO(3)$, although the results  of this article are much
more general  than that.  Since the action  is free,  a neighborhood of
the group orbit   $Gp_e$  may be  (with  respect to  the  left action)
equivariantly projected  to $SO(3)$ itself.  Let us realize $SO(3)$ as
the unit circle  bundle of the  2-sphere  in $\mathbb R^3$,  arranging
that the orbit $\exp(\xi t)p_e$ be projected to the the unit circle in
the tangent space of the vertical vector $\mathbf k$. The system, when
at the relative equilibrium  $p_e$,    and when viewed through    this
projection,  appears  merely as  a point,  say  $\mathcal P$, rotating
uniformly  on  the unit circle   in the tangent   space of $S^2$ at
$\mathbf k$,  and the other  relative equilibria $gp_e$, $g\in SO(3)$,
appear as   this  same motion but   reoriented  by  the  rotation $g$.
Suppose $p_e$  is formally  stable, and   suppose  the system is  then
perturbed from  the relative equilibrium $p_e$.   If $\mu_e\ne 0$ then
the  motion of the system  is  bound to remain near $\exp(\xi_e\mathbb
R)p_e$, as shown  in~\ct{PatrickGW-1992.1} (extensions  to the nonfree
case    may       be          found     in~\ct{OrtegaJPRatiuTS-1997.1}
and~\ct{LermanESingerSF-1997.1}).  The perturbation affects the motion
of the point $\mathcal P$ by superimposing on its original motion some
small vibration; however the   motion of $\mathcal  P$ does   not ever
carry that point far from its original circular path.

The situation is  different when $\mue=0$, for then  $G_{\mue}=SO(3)$,
and under perturbation the system is bound only to remain near $Gp_e$.
The  effect is that,   after a  perturbation   from $p_e$, the   point
$\mathcal P$ still moves nearly circularly in tangent spaces of $S^2$,
but the base point  of  the tangent  space slowly changes  location on
$S^2$.  In fact, as  shown  in~\ct{PatrickGW-1995.1}, the base   point
$\mathcal B$ moves as a  point charge under the  influence of a radial
magnetic   monopole.    This result     is  obtained  by  constructing
coordinates which make  an open neighborhood  of $Gp_e$ in $P$ into an
open neighborhood  of  $Z\bigl(T^*SO(3)\bigr)\times  \{\bar p_e\}$  in
$T^*SO(3)\times P_{\mu_e}$, where $P_{\mu_e}=J^{-1}(\mue)/G_{\mue}$ is
the  symplectic   reduced   space, $\bar p_e\in     P_{\mu_e}$  is the
equilibrium  of the   reduced    space corresponding  to   $p_e$,  and
$Z\bigl(T^*SO(3)\bigr)$ is the zero  section  of $T^*SO$.  In  the new
coordinates, to first order in  the momentum perturbation, the ``group
variables'' in $T^*SO(3)$ and the ``reduced variables'' in $P_{\mu_e}$
become decoupled; truncating  the higher order interaction terms gives
a new Hamiltonian  system   on $T^*SO(3)\times P_{\mu_e}$   called the
\defemph{drift system}.  The  drift   system has an  additional  $S^1$
\defemph{normal form symmetry}.  The  Hamiltonian system for  the slow
motion of the   base point $\mathcal B$ is   obtained from  the  drift
system  by   ignoring  the variables  on   $P_\mue$  and  reducing the
resulting Hamiltonian on $T^*SO(3)$ by  the normal form symmetry.  The
approximating,   truncated Hamiltonian of   the   drift system may  be
obtained merely by calculating the nilpotent part of the linearization
of the  vector  field $X_{H_\xie}(p)\equiv  X_H(p)-\xi_ep$  at  $p_e$,
where $\xi p\equiv (\xi)_P(p)$ denotes the infinitesimal generator of
$\xi_e$ at $p_e$.

Actually,  the  results   in~\ct{PatrickGW-1995.1}   depend    on    a
nonresonance condition that  the frequency $|\xi_e|$  cannot equal any
linearized frequency of the  reduced    system at $\bar  p_e$.    This
article considers the problem of constructing a Hamiltonian system for
the motion of  the base point $\mathcal B$  in the case where there is
such a \defemph{1:1 group-reduced resonance}.  The main result of this
article is that, in the presence of resonance, some of the reduced and
group  variables couple to first  order  in the momentum perturbation,
with    the     result  that $\mathcal  B$     moves     as  a charged
particle~\emph{with magnetic moment}  on a sphere under the  influence
of a  radial  magnetic monopole. However, this  ``particle'' regularly
exchanges its charge with its magnetic moment.

The construction of the Hamiltonian system which  models the motion of
$\mathcal  B$ in the  1:1 resonant case  generalizes and parallels the
nonresonant case.  In the course of the construction the linearization
of $X_{H_\xie}$  at $p_e$  is put into  a  certain  block  normal form
through   a Witt  or Moncreif    decomposition of  the  tangent  space
$T_{p_e}P$.  The  objective of this normal form   is the separation of
the  group and reduced  motions and it is not  the full normal form of
the linearization as  a infinitesimal symplectic operator.  The linear
normal form is   then extended to  a  neighborhood of  the group orbit
$Gp_e$ using  the the isotropic embedding/equivariant Darboux theorem.
These normal forms---for linearizations of relative equilibria and for
Hamiltonian systems near group orbits of relative equilibria---are the
second  main   contribution of  this article. In~\ct{PatrickGWRobertsRM-1999.1}
these normal forms have  been used to study the  structure of  the set  of
relative equilibria.

A new aspect that is absent in the nonresonant case but emerges in the
resonant case is the  appearance of a~\defemph{gauge group}.  There is
too much freedom inherent in the normal forms to exactly fit the needs
of first order agreement between the original and drift systems.  This
slack is taken up by the group $(\mathbb  R,+)$ and presents itself as
an        inherent    freedom  in     the    choice    of   the normal
forms.  This~\emph{gauge freedom} can  be  used to simplify the  drift
system, a motif which is quite useful in this work.

The  Hamiltonian system   modeling the   motion of   $\mathcal B$  is
approximate since  it is the result of  a  truncation, so there arises
the question of how well this model Hamiltonian reflects the behavior
of the  real  one.  To check  this  I  have numerically simulated  the
system of two (identical) axially symmetric rods which are joined by a
frictionless  ball-and-socket    joint,     but     otherwise     move
freely~(\ctt{PatrickGW-1989.1}{PatrickGW-1991.1}).       The  relative
equilibria of this system correspond to motions such that the two rods
spin  on their   axes   while otherwise maintaining constant    mutual
orientations, while the whole   assemblage  rotates about some   fixed
axis.   There  exist, due  to  the  system's multiple rotating  parts,
relative  equilibria which  are  formally stable  and  which have zero
(hence  nongeneric)    total-angular-momentum.  In  fact there    is a
continuum of relative equilibria with zero total-angular-momentum, and
it happens  that parameters may be chosen  so  that there are relative
equilibria at zero total-angular-momentum and at the 1:1 group-reduced
resonance. The third main  contribution of this  article is the actual
calculation of  the drift system for  the coupled rod  system near one
such   relative   equilibria,   and the  verification  of   the  drift
approximation by comparison  of  simulations of the coupled rod   system
itself and certain predictions of its drift system.

Here is   an overview of  this work.  I begin  in  Section~2.1 with an
analysis  of  the linearization~$dX_{H_\xie}(p_e)$  of $X_{H_\xie}$ at
$p_e$, by splitting  the  linearization into semisimple and  nilpotent
parts: $dX_{H_\xie}(p_e)= S_{p_e}+  N_{p_e}$.   The key  tool  in  the
analysis is a certain Moncrief or Witt decomposition  of~$T_{p_e}P$:
\begin{equation*}
T_{p_e} P\cong\Wred\oplus(\mathfrak g_\mue\oplus
\mathfrak g_\mue^*)\oplus T_\mue(G\mue).
\end{equation*} 
Whereas this decomposition is obtained in~\ct{PatrickGW-1995.1} for
the nonresonant case with the aid of the generalized eigenspaces of
$dX_{H_\xie}(p_e)$, and so was $dX_{H_\xie}(p_e)$-invariant, here the
decomposition is only $S_{p_e}$-invariant. Using this decomposition,
the nilpotent part~$ N_{p_e}$ acquires a block form which gives rise
to certain operators:
\begin{equation*}
N_{p_e}^{23}:\mathfrak g_\mue^*\rightarrow\mathfrak g_\mue,\quad 
N_{p_e}^{21}:T_{\bar p_e}P_\mue\rightarrow\mathfrak g_\mue,
\quad N_{p_e}^{13}:\mathfrak g_\mue^*\rightarrow T_{\bar p_e}P_\mue.
\end{equation*}
In the nonresonant case $N_{p_e}^{21}=0$ and $N_{p_e}^{13}=0$
necessarily. I expose, in Section~2.2, the special properties of
these operators, the most notable of which are that $N_{p_e}^{23}$ is
symmetric, $N_{p_e}^{13}$ and $N_{p_e}^{21}$ are dual, and that
these operators have certain commutation relations with operators
$\onm{ad}_\xie$ and $\onm{coad}_\xie$.

Next, in Section~3, I focus on the case where the reduced and group
spectra have intersection $\{\pm i\lambda_{p_e}\}$, a purely imaginary
eigenvalue and its conjugate, each of which occur with multiplicity~1
in both the reduced and group spectrum. The generalized eigenspace of
the reduced linearization corresponding to $\{\pm i\lambda_{p_e}\}$ is
symplectic and has dimension~2, and so is linearly symplectomorphic to
$\bigl(\mathbb R^2=\{x_1,x_2\},dx_1\wedge dx_2\bigr)$, and the
operator $N_{p_e}^{21}$ is zero on the sum of the complimentary
generalized eigenspaces. Consequently, the splitting of $T_{p_e}P$
refines and the operator $N_{p_e}^{21}$ may be replaced by another
operator $N_{p_e}^{211}:\mathbb R^2\rightarrow\mathfrak g_\mue$. After
this simplification, the splitting of $T_{p_e}P$ and the isotropic
embedding theorem~(\ct{MarsdenJE-1981.1},\ct{WeinsteinA-1997.1})
together give a map which transforms the Hamiltonian $H$, to first
order near $G_\mue p_e$, to the Hamiltonian
\begin{equation*}
\Hdrift(x,\alpha_g)\equiv\langle g^{-1}\alpha_g,\xi_e\rangle+
\frac12N_{p_e}^{23}(\alpha_g,\alpha_g)+
\langle g^{-1}\alpha_g, N_{p_e}^{211}x\rangle
\end{equation*}
near $\{0\}\times Z(T^*G_\mue)$ in the phase space $\mathbb R^2\times
T^*G_\mue$, where $ Z(T^*G_\mue)$ denotes the zero section. In
addition to the expected invariance under the left action of $G_\mue$,
$\Hdrift$ is invariant under a diagonal action of the toral subgroup
generated by $\xi_e$; this is the normal form symmetry.

In Section~4 I further assume that $G$ is the largest compact
continuous symmetry group for an ordinary mechanical system, namely
$SO(3)\times(S^1)^n$. After various manipulations, including an
Abelian reduction, the drift Hamiltonian becomes
\begin{equation*}
\Hdrift\equiv\frac12I_1\bigl({\pi_1}^2+{\pi_2}^2\bigr)+\frac12I_2{\pi_3}^2+
\kappa(\pi_1x_1+\pi_2x_2)+a\pi_3
\end{equation*}
on the phase space $\mathbb R^2\times
T^*SO(3)=\{(x,A,\pi)\}$. Here $I_1$ and $I_2$ come from
$N_{p_e}^{23}$, $\kappa$, a single coupling constant, is all that
remains of $N_{p_e}^{211}$, $a$ is a constant, and there are two
possibilities: the relative equilibrium can be either a ``+'' type or
a ``$-$'' type. In this context the gauge freedom has the effect of
making $I_1$ arbitrary and interest is focused on the dynamics near to
$x=0$, $\pi=0$. The Hamiltonian $\Hdrift$ is defined on a phase space
of dimension~$8$ and has symmetry $SO(3)\times S^1$, and so defines a
system which is completely integrable. Although a complete analysis
in the general case seems difficult, I analyze this system in
Section~4.1, where I show that the system has, for example, some
singular reduced phase spaces and spectrally unstable relative
equilibria with homoclinic connections. It is in Section~4.2 that I
show that the drift system reduced by the normal form symmetry can be
cast as a charged particle with magnetic moment moving on the sphere
while under the influence of a magnetic monopole.

Finally, in Section~5, I numerically investigate the dynamics of the
coupled rod system near one particular resonant relative equilibrium
and compare this dynamics with that of the drift system. I explore
three distinct regions of phase space: 1) within zero
total-angular-momentum, 2) near a stable relative equilibria of the
drift system, and 3) near a spectrally unstable relative equilibrium
at a singularity of the drift system. The comparison of the two
systems is hampered by the implicit nature of the coordinates relating
them, but in the first two comparisons agreement between the two
systems is obtained uneventfully. However, the singular points of
the third comparison involve an unexpected \emph{reconstuction phase
jump}; the situation is delicate and small perturbations are required
to elicit quantitative agreement between the two systems. Nevertheless
it becomes clear that many elements of the dynamics of the coupled rod
system near the resonant relative equilibrium are indeed captured by
the drift system.

\section{The linearization}
Here are the basic notations:
\begin{enumerate}
\item 
$p_e$ is a relative equilibrium,   $\xi_e$ is the generator of  $p_e$,
and the momentum of $p_e$ is $\mue\equiv J(p_e)$.
\item 
$G_\mue$ is the isotropy group of $\mue$ under the coadjoint action
of $G$ on $\mathfrak g^*$, and $\mathfrak g_\mue$ is the Lie algebra
of $G_\mue$. $\onm{CoAd_g}\equiv(\onm{Ad}_{g^{-1}})^*$, $g\in G$, and
$\onm{coad}_\xi\equiv-(\onm{ad}_\xi)^*$, $\xi\in\mathfrak g$.
\item 
$H_{\xi_e}\equiv H-J_{\xi_e}$, so that $H_{\xi_e}$ has a critical
point at $p_e$. The Hessian of $\Hxie$ at $p_e$ will be denoted by
$d^2\Hxie(p_e)$ and the linearization of $X_\Hxie$ at $p_e$ will be
denoted $dX_\Hxie(p_e)$.
\item
Without loss of generality, since everything is local near to $p_e$ and
$p_e$ is regular,
the Marsden-Weinstein symplectic reduced space $(P_\mue,\omega_\mue)$
exists; let $\pi_\mue:J^{-1}(\mue)\rightarrow P_\mue$ be the
projection. The reduced system has this phase space with Hamiltonian
$H_\mue$ defined by $H_{\mu_e}\pi_{\mu_e}=H|J^{-1}(\mue)$. Also,
$\bar p_e\equiv\pi_{\mu_e}(p_e)$, and the linearization of the reduced
system at $\bar p_e$ is denoted by $dX_{H_{\mu_e}}(\bar p_e)$. Suppose that
$p_e$ is regular, which means $\xi  p_e\neq0$ for all $\xi\in\mathfrak
g$.    
\item   The    \defemph{reduced  spectrum} is the    spectrum   of the
linearization $dX_{H_{\mu_e}}(\bar p_e)$ at  $\bar p_e$ of the reduced
system.   The   \defemph{group   spectrum} is  the   the   spectrum of
$\onm{ad}_\xie:\mathfrak g\rightarrow\mathfrak g$).
\end{enumerate}
Let $dX_\Hxie(p_e)= N_{p_e}+ S_{p_e}$ be  the Jordan decomposition  of
$dX_\Hxie(p_e)$ into its semisimple part  $S_{p_e}$ and nilpotent part
$ N_{p_e}$. The aim of this section is an analysis of $dX_\Hxie(p_e)$,
focusing  on its nilpotent  part $N_{p_e}$. \emph{Do not} assume that
the reduced spectrum and the group spectrum are disjoint.

\subsection{The splitting of $T_{p_e}P$}
I begin by deriving a ``Moncreif'' 
splitting~(\ct{MarsdenJE-1981.1}) of $T_{p_e}P$ which is
slightly weaker than its analogue in~\ct{PatrickGW-1995.1}, in that the
``reduced'' part (below $\Wred$) of the splitting is not (and cannot,
in general) be constructed to be $dX_\Hxie(p_e)$ invariant. The
details are similar to those in~\ct{PatrickGW-1995.1}; I will not
belabor them here. The subspace $\onm{ker}dJ(p_e)$ is
$S_{p_e}$-invariant, so one can choose an $S_{p_e}$-invariant subspace
$\Wred$ such that
\begin{equation*}\onm{ker}dJ(p_e)=\Wred\oplus\mathfrak g_\mue p_e.
\end{equation*}
The subspace $\Wred$ is symplectic, and
$T\pi_\mue|\Wred:\Wred\rightarrow T_{\bar p_e}P_\mue$ is a linear
symplectomorphism. Choose an $\onm{Ad}$-invariant complement
$\mathfrak b$ to $\mathfrak g_\mue$, so that
\begin{equation*}
\mathfrak g=\mathfrak b\oplus\mathfrak g_\mue.
\end{equation*}
The subspace $\mathfrak bp_e$ is symplectic, and $\mathfrak
g_\mue\subseteq(\Wred\oplus\mathfrak bp_e)^{\omega\perp}$ as an
$S_{p_e}$-invariant \emph{Lagrangian} subspace. Choose an
$S_{p_e}$-invariant Lagrangian complement\footnote{Since $S_{p_e}$ is
semisimple, every $S_{p_e}$-invariant Lagrangian subspace has an
$S_{p_e}$-invariant Lagrangian complement. This is a general fact;
the proof is a simple modification of the argument found at the bottom
of page~401 of~\ct{PatrickGW-1995.1}.} $Z$ to $\mathfrak g_\mue$ in
$(\Wred\oplus\mathfrak bp_e)^{\omega\perp}$, giving the
$S_{p_e}$-invariant splitting
\begin{equation}\lb{1}
T_{p_e}P=\Wred\oplus\mathfrak g_\mue p_e\oplus Z\oplus\mathfrak bp_e.
\end{equation}
As already stated, $T\pi_\mue$ is a symplectomorphism between $\Wred$
and $T_{p_e}P_\mue$. Also, $dJ(p_e)$ is a linear isomorphism between
$Z$ and $\mathfrak g_\mue^*$, and a symplectomorphism between
$\mathfrak b p_e$ and the tangent space $T_\mue(G\mue)$ at $\mue$ to
the coadjoint orbit $G\mue$. Thus~\rf{1} becomes, through these
identifications, the $ S_{p_e}$-invariant splitting,
\begin{equation}\lb{2}
T_{p_e}P=T_{\bar p_e}P_\mue\oplus\mathfrak g_\mue\oplus \mathfrak g_\mue^*
\oplus T_\mue(G\mue).
\end{equation}

With respect to the decomposition~\rf{2}, let $dX_\Hxie(p_e)$,
$S_{p_e}$, and $N_{p_e}$ have blocks $[A^{ij}]$, $[S^{ij}]$, and
$[N_{p_e}^{ij}]$, respectively; the block form of $\omega(p_e)$ becomes
\begin{equation}\lb{3}
\omega(p_e)=\left[\begin{array}{cccc}
\omega_\mue(\bar p_e)&0&0&0\\
0&0&\Id&0\\
0&-\Id&0&0\\
0&0&0&\breve\omega_\mue(\mue)
\end{array}\right],\end{equation}
where $\omega_\mue$ is the reduced symplectic form of $P_\mue$ and
$\breve\omega_\mue$ is the Kostant-Souriau form on $G\mue$.
In~\rf{1}, the subspaces $\mathfrak g_\mue p_e$, $\mathfrak b p_e$,
and $\onm{ker}dJ(p_e)=\Wred\oplus\mathfrak g_\mue p_e$ are
$dX_\Hxie(p_e)$-invariant, and hence are $N_{p_e}$-invariant, and this
implies certain of the $A^{ij}$ and $N_{p_e}^{ij}$ vanish.
Using the identities (Proposition~5 of~\ct{PatrickGW-1995.1})
\begin{gather*}
X_\Hxie(p_e)=-\onm{coad}_{\xi_e}dJ(p_e)\\
dX_\Hxie(p_e)\eta(p_e)=-(\onm{ad}_{\xi_e}\eta)_P(p_e),\quad\eta\in
\frak g
\end{gather*}
to calculate the diagonal blocks of $dX_\Hxie(p_e)$, the Jordan
decomposition becomes
\begin{equation*}\begin{split}
dX_\Hxie(p_e)&=\left[\begin{array}{cccc}
dX_{H_\mue}(\bar p_e)&0&A^{13}&0\\
\rule{0pt}{11pt}A^{21}&-\onm{ad}_{\xi_e}&A^{23}&0\\
0&0&\onm{ad}^*_{\xi_e}&0\\
0&0&0&\onm{ad}^*_{\xi_e}
\end{array}\right]\\
&=
\left[\begin{array}{cccc}
S^{11}&0&0&0\\
0&S^{22}&0&0\\
0&0&S^{33}&0\\
0&0&0&S^{44}
\end{array}\right]
+{\renewcommand\arraystretch{1.15}\left[\begin{array}{cccc}
N_{p_e}^{11}&0&N_{p_e}^{13}&0\\
N_{p_e}^{21}&N_{p_e}^{22}&N_{p_e}^{23}&0\\
0&0&N_{p_e}^{33}&0\\
0&0&0&N_{p_e}^{44}
\end{array}\right]}.
\end{split}\end{equation*}
Since $S_{p_e}$ and $N_{p_e}$ commute, so do $S^{11}$ and $N_{p_e}^{11}$,
and $S^{11}$ is semisimple and $N_{p_e}^{11}$ is nilpotent since $S_{p_e}$
and $N_{p_e}$ are. Thus, $dX_{H_\mue}(\bar p_e)=S^{11}+N_{p_e}^{11}$ is the
Jordan decomposition of of the reduced linearization $dX_{H_\mue}(\bar
p_e)$, which is semisimple, since $p_e$ is formally
stable. Consequently, $S^{11}=dX_{H_\mue}(\bar p_e)$ and $N_{p_e}^{11}=0$.
Similarly, $N_{p_e}^{22}=N_{p_e}^{33}=N_{p_e}^{44}=0$, and
\begin{equation*}
S^{22}=-\onm{ad}_{\xi_e},\quad S^{33}=\onm{ad}^*_{\xi_e},
\quad S^{44}=\onm{ad}^*_{\xi_e}.\end{equation*}
Thus,
\begin{equation}\lb{4}
 S_{p_e}=\left[\begin{array}{cccc}
dX_{H_\mue}(\bar p_e)&0&0&0\\
\rule{0pt}{11pt}0&-\onm{ad}_{\xi_e}&0&0\\
0&0&(\onm{ad}_{\xi_e})^*&0\\
0&0&0&(\onm{ad}_{\xi_e})^*
\end{array}\right]
\end{equation}
and
\begin{equation}\lb{5}
 N_{p_e}={\renewcommand\arraystretch{1.15}\left[\begin{array}{cccc}
0&0&N_{p_e}^{13}&0\\
N_{p_e}^{21}&0&N_{p_e}^{23}&0\\
0&0&0&0\\
0&0&0&0
\end{array}\right]}.
\end{equation}
As will be seen immediately below, a large amount of information can
be discerned from structure of $N_{p_e}$ and $S_{p_e}$
visible in~\rf{4} and~\rf{5}.

\subsection{Properties of $N_{p_e}$ and its subblocks}
Directly from~\rf{5}, $( N_{p_e})^3=0$ so that $N_{p_e}$ has nilpotent
order at most~3. This upper bound on the nilpotent order occurs
because $dX_\Hxie(p_e)$ is semisimple when suitably restricted and
projected to the summands of~\rf{1}, so that $N_{p_e}$ must
``transport'' between those summands. Moreover, the ``direction'' of
this transport is ``one way'', since $\onm{ker}dJ(p_e)$ and $\mathfrak
g_\mue p_e$ are invariant subspaces, and the ``transport'' stops at
$\mathfrak g_\mue p_e$, on which $ N_{p_e}$ is zero:
\begin{gather}
\onm{Image} N_{p_e}\subseteq \onm{ker}dJ(p_e),\lb{6}\\
N_{p_e}\bigl(\onm{ker}dJ(p_e)\bigr)\subseteq\mathfrak g_\mue p_e,\nonumber\\
N_{p_e}|\mathfrak g_\mue p_e = 0.\lb{8}
\end{gather}
When iteratively acted upon by $N_{p_e}$, a vector can make
just~2~stops before annihilation: $\onm{ker}dJ(p_e)$ and $\mathfrak
g_\mue p_e$; hence the nilpotent order of $N_{p_e}$ is a most~3.

By~\rf{8}, 
\begin{equation*}
\onm{ker}T_{p_e}\pi_\mue=\mathfrak g_\mue p_e\subseteq\onm{ker} \bigl(
N_{p_e}|\onm{ker}dJ(p_e)\bigr),
\end{equation*} 
so $N_{p_e}^{21}:T_{\bar p_e}P_\mue\rightarrow\mathfrak g$ is the
unique linear map such that
\begin{equation}\lb{24}\begin{diagram}
\node{\onm{ker}dJ(p_e)}\arrow{e,t}{N_{p_e}}\arrow{s,l}{T_{p_e}\!\pi_\mue}
 \node{\mathfrak g_\mue p_e}\arrow{e,t}{\xi p_e\mapsto\xi}
 \node{\mathfrak g_\mue}\\
\node{T_{\bar p_e}P_\mue}\arrow{ene,b}{N_{p_e}^{21}}
\end{diagram}\end{equation}
Since this diagram does not depend on choices made in the construction
of the splitting~\rf{1}, neither does $N^{21}$. Also, $N_{p_e}^{21}$
may be calculated merely by reference to the diagram; it is not
necessary to calculate all summands of~\rf{1}. Similarly, by~\rf{6},
and since, by~\rf{5},
\begin{equation*}\onm{ker}dJ(p_e)\subseteq\onm{ker}(T\pi_\mue N_{p_e}),
\end{equation*}
 $N_{p_e}^{13}:\mathfrak g_\mue^*\rightarrow
T_{\bar p_e}P_\mue$ is the unique linear map such that 
\begin{equation*}\begin{diagram}
\node{T_{p_e}P}\arrow{e,t}{dJ(p_e)}\arrow{s,l}{N_{p_e}}
 \node{\frak g^*}\arrow{e}\node{\mathfrak g_\mue^*}
 \arrow[2]{sw,b}{N_{p_e}^{13}}\\
\node{\onm{ker}dJ(p_e)}\arrow{s,l}{T_{p_e}\!\pi_\mue}\\
\node{T_{\bar p_e}P_\mue}
\end{diagram}\end{equation*}

Since $dX_\Hxie(p_e)$ is infinitesimally symplectic, so is $N_{p_e}$,
so that
\begin{equation}\lb{9}( N_{p_e})^t\omega(p_e)+\omega(p_e)N_{p_e}=0,
\end{equation} 
and inserting ~\rf{3} and~\rf{5} into~\rf{9} gives
\begin{equation}\lb{10}N_{p_e}^{21}=-(N_{p_e}^{13})^*\omega_\mue^\flat,
\qquad(N_{p_e}^{23})^*=N_{p_e}^{23},
\end{equation} 
so $N_{p_e}^{23}$ is symmetric, and in a certain sense, $N_{p_e}^{13}$
and $N_{p_e}^{21}$ are dual. Similarly, using~\rf{4} and~\rf{5} to
write out what it means for $S_{p_e}$ and~$ N_{p_e}$ to commute,
\begin{gather}
\lb{14}dX_{H_\mue}(\bar p_e)N_{p_e}^{13}=-N_{p_e}^{13}
 \onm{coad_{\xi_e}},\\
\lb{15}N_{p_e}^{21}dX_{H_\mue}(\bar p_e)=-\onm{ad}_{\xi_e}N_{p_e}^{21},\\
\lb{16}\onm{ad}_{\xi_e}N_{p_e}^{23}=N_{p_e}^{23}\onm{coad}_{\xi_e}.
\end{gather}
These properties will yield the important \emph{normal form symmetry}
of the drift system.

Temporarily let $\onm{pr}_2$ be the projection onto the second factor
of~$\onm{ker}dJ(p_e)=\Wred\oplus \mathfrak g_\mue p_e$. Then
\begin{equation*}
\onm{ker}\bigl(dJ(p_e)|(\Wred+\mathfrak bp_e)^{\omega\perp}\bigr)
=\mathfrak g_\mue p_e\subseteq\onm{ker}\bigl(\onm{pr}_2N_{p_e}|
(\Wred+\mathfrak bp_e)^{\omega\perp}\bigr),
\end{equation*}
so $N_{p_e}^{23}:\mathfrak g_\mue\times\mathfrak
g_\mue\rightarrow\mathbb R$ is the unique bilinear from such that
\begin{equation}\lb{17}\begin{diagram}
\node{(\Wred+\mathfrak bp_e)^{\omega\perp}}\arrow{e,t}{dJ(p_e)}
 \arrow{s,l}{N_{p_e}}
 \node{\mathfrak g_\mue^*}\arrow{ssw,b}{(N_{p_e}^{23})^\flat}\\
\node{\onm{ker}dJ(p_e)}\arrow{s,l}{\onm{pr}_2}\\
\node{\mathfrak g_\mue}
\end{diagram}\end{equation}
Thus, $N^{23}$ may be calculated
without calculating the Lagrangian complement $Z$ of~\rf{1}.

However, $N_{p_e}^{23}$ \emph{does} depend on the particular choice or
$\Wred$, even if $\Wred$ is an $S_{p_e}$-invariant complement of
$\mathfrak g_\mue p_e$ in $\onm{ker}dJ(p_e)$. To see this dependence,
choose a linear map $A:\Wred\rightarrow \mathfrak g_\mue$, which
commutes with $S_{p_e}$, and then choose a new $\Wred$ where
\begin{equation}\lb{12}
\Wred^{\,\prime}=\bset{w+(Aw)p_e}{w\in\Wred}.
\end{equation}
Also, choose a map $B:\mathfrak b\rightarrow\mathfrak g_\mue$ which
intertwines the adjoint action and then choose a new $\mathfrak
b^\prime$ by
\begin{equation}\lb{11}
\mathfrak b^\prime=\bset{\eta+B\eta}{\eta\in\mathfrak b}.
\end{equation}
Given $\nu\in\mathfrak g_\mue^*$, the process of calculating
$({N_{p_e}^{23}}^\prime)^\flat\nu$, given by the analog of~\rf{17} for
the new choices, is
\begin{enumerate}
\item 
pick $z^\prime\in (\Wred^{\,\prime}+\mathfrak b^\prime
p_e)^{\omega\perp}$ such that $dJ(p_e)z^\prime=\nu$;
\item calculate $N_{p_e}z^\prime$;
\item 
project $N_{p_e}z^\prime$ to $\mathfrak g_\mue p_e\equiv\mathfrak
g_\mue$ using the splitting of $\onm{ker}dJ(p_e)$ defined by
$\Wred^{\,\prime}$, with result $({N_{p_e}^{23}}^\prime)^\flat\nu$.
\end{enumerate}
It is straightforward from~\rf{3},~\rf{12}, and~\rf{11}, that
\begin{equation*}
(\Wred^{\,\prime}+\mathfrak b^\prime p_e)^{\omega\perp}
=\bset{\omega_\mue^\sharp A^*\mu\oplus\xi\oplus\mu\oplus\breve\omega_\mue
^\sharp B^*\mu}{\xi\in\mathfrak g_\mue,\mu\in\mathfrak g_\mue^*},
\end{equation*}
so given $\nu\in\mathfrak g_\mue^*$, let
\begin{equation*}
z^\prime=\omega_\mue^\sharp A^*\nu\oplus0\oplus\nu\oplus\breve
\omega_\mue^\sharp B^*\nu.\end{equation*}
Then 
\begin{equation*}\begin{split}
 N_{p_e} z^\prime&=
N_{p_e}^{13}\nu\oplus(N_{p_e}^{21}\omega_\mue^\sharp A^*\nu
+(N_{p_e}^{23})^\flat\nu)\oplus0\oplus0\\
&=N_{p_e}^{13}\nu\oplus AN_{p_e}^{13}\nu\oplus0\oplus0\\
&\qquad\mbox{}+0\oplus(N_{p_e}^{21}\omega_\mue^\sharp A^*\nu
-AN_{p_e}^{13}\nu+(N_{p_e}^{23})^\flat\nu)\oplus0\oplus0,
\end{split}\end{equation*}
which by~\rf{12} is the appropriate decomposition of
$N_{p_e}z^\prime$. Thus, using~\rf{10},
\begin{equation}\lb{13}\begin{split}
({N_{p_e}^{23}}^\prime)^\flat\nu&=
N_{p_e}^{21}\omega_\mue^\sharp A^*\nu
-AN_{p_e}^{13}\nu+(N_{p_e}^{23})^\flat\nu\\
&=-\bigl((N_{p_e}^{13})^*A^*+AN_{p_e}^{13}\bigr)\nu
+(N_{p_e}^{23})^\flat\nu.
\end{split}\end{equation}
This freedom to adjust $N_{p_e}^{23}$ by adjusting $\Wred$ will give
the important \emph{gauge freedom} of the drift
system.

This gauge freedom can be encoded as a group, as follows. Temporarily
set $\mathbb E_{p_e}=T_{\bar p_e}P_\mue\times\frak g_\mue\times\frak
g_\mue^*\times T_\mue G_\mue$, and fix one particular splitting of
type~\rf{1}, giving a symplectomorphism, say
$\phi_0:T_{p_e}P\rightarrow\mathbb E_{p_e}$. Given
$A:\Wred\rightarrow\frak g_\mue$ and $B:\frak b\rightarrow\frak
g_\mue$ as above, there is another symplectomorphism
$\phi_{A,B}:T_{p_e}P\rightarrow\mathbb E_{p_e}$, and hence a unique
$\Delta_{A,B}$ in the symplectic group $Sp(\mathbb E_{p_e})$ such that
\begin{equation*}\begin{diagram}
\node[2]{T_{p_e}P}\arrow{sw,l}{\phi_0}\arrow{se,t}{\phi_{A,B}}\\
\node{\mathbb E_{p_e}}\arrow[2]{e,b}{\Delta_{A,B}}\node[2]{\mathbb E_{p_e}}
\end{diagram}\end{equation*}
Conversely, $A$ is determined by $\Delta_{A,B}$, since if $x\in
T_{\bar p_e}P_\mue$, $\xi\in\frak g_\mue$, and $w\in\Wred$ is such
that $\phi_0(w)=x$, then
\begin{multline*}
\Delta_{A,b}(x)=\phi_{A,B}(w)
=\phi_{A,B}\bigl(w+(Aw)p_e-(Aw)p_e\bigr)\\=
T\pi_\mue\bigl(w+(Aw)p_e\bigr)\oplus(-Aw)
=x\oplus(-Aw).
\end{multline*}
Similarly $B$ is determined by $\Delta_{A,B}$, so the pairs $(A,B)$
have a natural group structure given by the injection
$(A,B)\mapsto\phi_{A,B}$, and a simple calculation gives
\begin{equation*}
\phi_{A+A^\prime,B+B^\prime}=\phi_{A,B}\circ\phi_{A^\prime,B^\prime}.
\end{equation*}
Thus, at this level, the \emph{gauge group} is the additive Abelian
group of pairs $(A,B)$ such that $A$ commutes with $S_{p_e}$ and $B$
commutes with $\onm{ad}_{\xi_e}$.

\section{The simplest resonant case}
All the above has been established without any particular presumption
on the way that the group and reduced spectra intersect. There are
some apriori features of the spectra: the spectra are purely imaginary
(or zero) and are invariant under change of sign. Since $p_e$ is
presumed to be formally stable, zero cannot occur in the reduced
spectrum, and hence not in the intersection of the spectra. Thus, the
following \defemph{resonance condition} is the simplest possible case
beyond an empty intersection:
\begin{equation*}\parbox{.8\displaywidth}{
\em The intersection of the spectrum of $\onm{ad}_\xie$ and the
spectrum of $dX_{H_\mue}(\bar p_e)$ is $\{\pm i\lambda_{p_e}\}$,
$\lambda_{p_e}>0$, and each of the eigenvalues $\pm i\lambda_{p_e}$
occur in each of these operators with multiplicity one.}
\end{equation*}
While many of the results below remain true or have analogues in a
more relaxed environment, for simplicity, this resonance condition
will be assumed for the remainder of this article.

\subsection{The linearization}
Given the resonance assumption, it is natural to consider certain
spectral splittings of the factors $\Wred$ and~$\mathfrak g_\mue$
in~\rf{2}. In particular, let $\Wred^1$ be the generalized eigenspace
of $\pm i\lambda_{p_e}$ for $dX_{H_\mue}(\bar p_e)$, and $\Wred^0$ be
the sum of the other generalized eigenspaces, so that $\Wred^1$ is two
dimensional, symplectic, and there is the symplectic splitting
\begin{equation}\lb{101}
T_{\bar p_e}P_\mue=\Wred^0\oplus\Wred^1.
\end{equation}
Similarly, let $\mathfrak g_\mue^1$ be the generalized eigenspace of $\pm
i\lambda_{p_e}$ for the linear map $\onm{ad}_\xie$, and $\mathfrak
g_\mue^0$ be the sum of the other generalized eigenspaces, so $\mathfrak
g_\mue^1$ is two dimensional, and 
\begin{equation}\lb{102}
\mathfrak g_\mue=\mathfrak g_\mue^0\oplus\mathfrak g_\mue^1.
\end{equation}
Of course, there is also then the dual
splitting $\mathfrak g_\mue^*=\mathfrak g_\mue^{0*}\oplus\mathfrak g_\mue^{1*}$.

So we have a refinement of the splitting~\rf{2} into
``nonresonant'' (with superscript~0) and ``resonant'' (superscript~1,
dimension~2) parts:
\begin{align}
T_{p_e}P&=\Wred^0\oplus\Wred^1\oplus\mathfrak g_\mue\oplus
\mathfrak g_\mue^*\oplus T_\mue(G\mue)\lb{27}\\
&=\Wred^0\oplus\Wred^1\oplus\mathfrak g_\mue^0\oplus\mathfrak g_\mue^1\oplus
\mathfrak g_\mue^{0*}\oplus\mathfrak g_\mue^{1*}\oplus T_\mue(G\mue).\nonumber
\end{align}
As one might expect, $N_{p_e}^{13}$ and $N_{p_e}^{21}$
localize to the resonant parts:
\begin{alignat}{2}
&\Wred^0\subseteq\onm{ker}N_{p_e}^{21},&
\qquad &\onm{Image}N_{p_e}^{21}\subseteq\mathfrak g_\mue^1,\lb{19}\\
&\mathfrak g_\mue^0\subseteq\onm{ker}N_{p_e}^{13},&\qquad&\onm{Image} 
N_{p_e}^{13}\subseteq\Wred^1.\lb{20}
\end{alignat}

Indeed, to show the first of~\rf{19}, it suffices to find a subspace,
say ${\Wred^0}^{\!\!\!\prime}\;$ of $\onm{ker}dJ(p_e)$ such that
\begin{equation}\lb{21}
T_{p_e}\pi_\mue{\Wred^0}^{\!\!\!\prime}\;=\Wred^0\quad\mbox{and}\quad
 N_{p_e}|{\Wred^0}^{\!\!\!\prime}\; = 0.
\end{equation}
The result then follows by using ${\Wred^0}^{\!\!\!\prime}\;$ to
reverse the vertical arrow of~\rf{24}. The first of~\rf{21} can be
assured by setting ${\Wred^0}^{\!\!\!\prime}\;$ to be the sum of the
of nonresonant (i.e. not $\pm i\lambda_{p_e}$) generalized
eigenspaces of $dX_\Hxie(p_e)$, and the second of~\rf{21} follows
since ${\Wred^0}^{\!\!\!\prime}\;$ is $dX_\Hxie(p_e)$-invariant and
that $dX_\Hxie(p_e)|{\Wred^0}^{\!\!\!\prime}\;$ is semisimple. For
the second of~\rf{19}, $N_{p_e}$ maps the $\pm i\lambda_{p_e}$
generalized eigenspace of $dX_\Hxie(p_e)$ to itself, since
$dX_\Hxie(p_e)$ does that. Thus the image of $N_{p_e}$ is
contained in the intersection of that generalized eigenspace with
$\mathfrak g_\mue p_e$, which is exactly $\mathfrak g_\mue^1 p_e$.
From~\rf{19} one gets~\rf{20} by using the duality~\rf{10} between
$N_{p_e}^{21}$ and~$N_{p_e}^{13}$.

With respect to~\rf{27}, and in view of~\rf{5},~\rf{19}
and~\rf{20}, the linear map $N_{p_e}$ has the form
\begin{equation*}
 N_{p_e}=\left[\begin{array}{ccccc}
0&0&0&0&0\\
0&0&0&N_{p_e}^{131}&0\\
0&N_{p_e}^{211}&0&N_{p_e}^{23}&0\\
0&0&0&0&0\\
0&0&0&0&0\end{array}\right]
\end{equation*}
where $N_{p_e}^{211}\equiv N_{p_e}^{21}|\Wred^1$ and where
$N_{p_e}^{131}$ is $N_{p_e}^{13}$ regarded as a map into $\Wred^1$.
Let the quadratic Hamiltonian for 
the infinitesimally linear map
$ N_{p_e}$ be $\Hnil$. With respect to~\rf{27}, the array for $\Hnil$ is
\begin{multline*}
\left[\begin{array}{ccccc}
0&0&0&0&0\\
0&0&0&N_{p_e}^{131}&0\\
0&N_{p_e}^{211}&0&N_{p_e}^{23}&0\\
0&0&0&0&0\\
0&0&0&0&0\end{array}\right]^*
\left[\begin{array}{ccccc}
\omega_\mue^0(\bar p_e)&0&0&0&0\\
0&\omega_\mue^1(\bar p_e) &0&0&0\\
0&0&0&\Id&0\\
0&0&-\Id&0&0\\
0&0&0&0&\breve\omega_\mue(p_e)
\end{array}\right]\\
=\left[\begin{array}{ccccc}
0&0&0&0&0\\
0&0&0&(N_{p_e}^{211})^*&0\\
0&0&0&0&0\\
0&(N_{p_e}^{131})^*\omega_\mue^1(\bar p_e)&0&N_{p_e}^{23}&0\\
0&0&0&0&0\end{array}\right],\qquad\qquad
\end{multline*}
where $\omega_\mue^i(\bar p_e)$ denotes the symplectic form on
$\Wred^i$, $i=1,2$. Using~\rf{16}, and the fact that the projections
from the splitting $\mathfrak g_\mue^*=\mathfrak g_\mue^{0*}\oplus\mathfrak
g_\mue^{1*}$ are certain polynomials in $\onm{coad}_\xie$, the bilinear
form $N_{p_e}^{23}$ block diagonalizes over that splitting. Denote
the blocks by $N_{p_e}^{230}$ and $N_{p_e}^{231}$. Then, using the
variables
\begin{equation*}
(x,\mu)=(x,\mu^0,\mu^1)\in\Wred^1\times\mathfrak g_\mue^*=\Wred^1\times\mathfrak g_\mue^{0*}\times\mathfrak g_\mue^{1*},
\end{equation*}
we have
\begin{equation}\lb{33}\begin{split}
\Hnil&=\frac12N_{p_e}^{23}(\mu,\mu)+
\langle\mu,N_{p_e}^{211}x\rangle\\
&=\frac12N_{p_e}^{230}(\mu^0,\mu^0)+\frac12N_{p_e}^{231}(\mu^1,\mu^1)
+\langle\mu^1,N_{p_e}^{211}x\rangle.
\end{split}\end{equation}

I now construct the normal form symmetry for $\Hnil$. Let $T_\xie$ be
the closure of $\exp(\mathbb R\xi_e)$, so $T_\xie$ is Abelian, is
generically a maximal torus of $G_\mue$, and $T_\xie$ naturally acts
on both $\mathfrak g_\mue$ and $\mathfrak g_\mue^*$, with invariant
subspaces $\mathfrak g_\mue^i$ and $\mathfrak g_\mue^{i*}$, $i=1,2$.
By~\rf{16}, $N_{p_e}^{23}$ is $T_{\xi_e}$-invariant, and so~\rf{33} is
invariant under $T_{\xi_e}$ restricted to the $\mu^0$ variables alone.
After this, it is the last two terms of~\rf{33} that will dictate its
further symmetries. Now $S^1$ also acts\footnote{Define
$\theta^\wedge=\theta\left(\begin{array}{cc}0&-1\\1&0\end{array}\right)$
so that $S^1=SO(2)$ and $\exp$ is the usual matrix exponential.} on
$\mathfrak g_\mue^1$ by
\begin{equation*}
\exp(\theta^\wedge)\eta\equiv\exp\left(\frac\theta{\lambda_{p_e}}\xie
\right)\eta,\quad\theta\in\mathbb R,
\end{equation*}
which is merely a renormalization of the action of $T_{\xi_e}$ on
$\mathfrak g_\mue^1$, and by duality there is a corresponding action
of $S^1$ on $\mathfrak g_\mue^*$. Also there is the following
symplectic action of $S^1$ on $\Wred^1$:
\begin{equation}\lb{107}
\exp(\theta^\wedge) x\equiv\exp\left(\frac\theta{\lambda_{p_e}}
dX_{H_\mue}(\bar p_e)\right)x,\quad\theta\in\mathbb R.
\end{equation} 
By~\rf{14} and~\rf{15}, $N_{p_e}^{131}$ and $N_{p_e}^{211}$
reverse-intertwine the two $S^1$ actions, so that
\begin{equation*}
\langle{\exp(\theta^\wedge)}^{-1}\mu^1, N_{p_e}^{211}\exp(\theta^\wedge)
x\rangle =\langle{\exp(\theta^\wedge)}^{-1}\mu^1,
 {\exp(\theta^\wedge)}^{-1} N_{p_e}^{211}x\rangle=\langle\mu^1,x\rangle.
\end{equation*}
Consequently, the \defemph{normal form symmetry} will be
$T_{\xi_e}\times S^1$ where the first factor acts on the $\mu^0$
variables via the coadjoint action and the second factor acts by
\begin{equation*}
a\cdot(x,\mu^1)\equiv(a^{-1}x,a\mu),\quad a=\exp(\theta^\wedge).
\end{equation*}

The normal form symmetry restricts the quadratic Hamiltonian~$\Hnil$,
as follows. There is a basis of the linear space
$\Wred^1=\{(x_1,x_2)\in\mathbb R^2\}$ so that the symplectic
form is $dx^1\wedge dx^2$ and the $S^1$ action on $\Wred^1$ is either
$\bigl(\exp(\theta^\wedge),x\bigr)\mapsto \exp(\theta^\wedge)x$ or
$\bigl(\exp(\theta^\wedge),x\bigr)\mapsto \exp(-\theta^\wedge)x$ (i.e.
counterclockwise or clockwise rotation). Using the
$\onm{Ad}$-invariant metric of $\mathfrak g$, choose a orthonormal
basis for the space $\mathfrak g_\mue^1$, so that $\mathfrak
g_\mue^1=\{(\pi_1,\pi_2)\in\mathbb R^2\}$ and the $S^1$
action on $\mathfrak g_\mue^1$ is a rotation in the opposite sense as
the $S^1$ action on $\Wred$. The map $N_{p_e}^{211}$ is highly
restricted: since it reverse-intertwines these two reverse $S^1$
actions, there is a $\kappa$ such that
\begin{equation}
N_{p_e}^{211}=\kappa\left[\begin{array}{cc}1&0\\0&1\end{array}\right]
\qquad\mbox{or}\qquad
N_{p_e}^{211}=\kappa\left[\begin{array}{cc}0&-1\\1&0\end{array}\right].
\lb{40}
\end{equation} 
Equation~\rf{40} defines $\kappa$, so that the map $N^{211}$ has the
effect of giving rise to this single ``coupling constant'' and one of
the two possible choices in~\rf{40}. Without loss of generality,
by the symplectic isomorphisms 
$(x_1,x_2)\mapsto(-x_1,-x_2)$ and
$(x_1,x_2)\mapsto(-x_2,x_1)$ 
on $\Wred^1$ ,
$\kappa\ge0$
and the second choice in~\rf{40} may be discarded, respectively.
Then, since $N^{23}_{p_e}$ is
$S^1$-symmetric, a constant $I_1$ may be obtained from $N_{p_e}^{231}$
such that
\begin{equation}\lb{100}
\Hnil=
\frac12N_{p_e}^{230}(\mu^0,\mu^0)
+\frac12I_1({\pi_1}^2+{\pi_2}^2)
+\kappa(\pi_1x_1+\pi_2x_2).
\end{equation}

According to the above conventions, $N^{211}$ reverse intertwines
reverse-sense $S^1$ actions on the $(x_1,x_2)$ variables and the
$(\pi_1,\pi_2)$ variables, and the second action has the same sense as
the action of $\onm{ad}_\xie$. Consequently, {\em there are two
distinct possibilities:} the $S^1$ action on the $(x_1,x_2)$ variables
can be counterclockwise or clockwise. To avoid unnecessary signs that
would otherwise appear later, I will redefine the $S^1$ actions in the
following way: if the $S^1$ action on the $(x_1,x_2)$ variables is
clockwise then reverse it, and both actions are then counterclockwise
and the $S^1$ action on the $(\pi_1,\pi_2)$ variables is in the same
sense as the action of $\onm{ad}_\xie$. If on the other hand the $S^1$
action on the $(x_1,x_2)$ variables is counterclockwise then the
action of the $(\pi_1,\pi_2)$ variables is clockwise, and is to be
reversed, so that it has the opposite sense to the action of
$\onm{ad}_\xie$. Thus, by these conventions, both actions are
counterclockwise, $N^{211}$ is intertwining, and the distinctness of
the two cases appears as one of two possibilities: the $S^1$ action on
the $(\pi_1,\pi_2)$ variables might be in the same sense as the action
of $\onm{ad}_\xie$ (I call this the ``$+$'' case, and it occurs when
the original action on the $(x_1,x_2)$ variables is clockwise) or the
opposite sense as the action of $\onm{ad}_\xie$ (the ``$-$'' case,
occurring when the original action on the $(x_1,x_2)$ variables is
counterclockwise). Although it may appear that same-sense vs.
opposite-sense with respect to $\onm{ad_\xie}$ holds some
significance, similar to, for example, the difference between a $1:1$
and $1:-1$ resonance as in~\ctt{KummerM-1976.1}{KummerM-1978.1}, this
is not so: changing the the symplectic form $dx_1\wedge dx_2$ to
$dx_2\wedge dx_1$ will exchange my ``$+$'' and ``$-$'' cases. These
two cases are qualitatively identical and the subsequent theory is
ambidextrous with respect to them.

Now I will determine the effect on~\rf{100} of choosing various
subspaces $\Wred$, or equivalently, various operators~$A$ commuting
with $S_{p_e}$ in~\rf{12}. With respect to the
decompositions~\rf{101} and~\rf{102} let $A$ have the block form
\begin{equation*}
A=\left[\begin{array}{cc}A^{11}&A^{12}\\A^{21}&A^{22}\end{array}\right].
\end{equation*} 
To determine which $A$ commute with $S_{p_e}$, note that $S_{p_e}$ has
a block diagonal form on $\Wred^0\oplus\Wred^1\oplus\mathfrak
g_\mue^0\oplus\mathfrak g_\mue^1$ since~\rf{101} and~\rf{102} are
generalized eigenspace decompositions for $S_{p_e}$; let the blocks be
$S^{110}$, $S^{111}$, $S^{220}$, and $S^{221}$, so $AS=SA$ becomes
\begin{equation}\lb{105}
\left[\begin{array}{cc}A^{11}&A^{12}\\A^{21}&A^{22}\end{array}\right]
\left[\begin{array}{cc}S^{110}&0\\0&S^{111}\end{array}\right]=
\left[\begin{array}{cc}S^{220}&0\\0&S^{221}\end{array}\right]
\left[\begin{array}{cc}A^{11}&A^{12}\\A^{21}&A^{22}\end{array}\right].
\end{equation}
When multiplied, the top left blocks of~\rf{105} give
$A^{11}S^{110}=S^{220}A^{11}$, which is possible only if $A^{11}$ is
zero, since $S^{110}$ and $S^{220}$ have no common spectrum, and
similarly $A^{12}=0$ and $A^{21}=0$. However, $S^{111}$ and $S^{221}$
both have spectrum $\{\pm\lambda_e\}$, so one cannot conclude
$A^{22}=0$, but only that $A^{22}S^{111}=S^{221}A^{22}$. Since
$S^{111}$ is exactly the linearization~$dX_{H_\mue}(\bar p_e)$
restricted to $\Wred^1$, which by~\rf{107} defined the $S^1$ action on
$\Wred^1$, which is counterclockwise or clockwise, we have
\begin{equation*}
S^{111}=\pm\sqrt{\lambda_e}\left[\begin{array}{cc}0&1\\-1&0
\end{array}\right].
\end{equation*}
Since $S^{221}$ is $-\onm{ad}_{\xi_e}$ on $\mathfrak g_\mue^1$, and
$\onm{ad}_{\xi_e}$ generates the $S^1$ action of opposite sense to
that generated by $S^{111}$, we have
\begin{equation*}
S^{221}=\pm\sqrt{\lambda_e}\left[\begin{array}{cc}0&1\\-1&0
\end{array}\right].
\end{equation*}
By letting
\begin{equation*}
A^{22}=\left[\begin{array}{cc}a^{11}&a^{12}\\a^{21}&a^{22}
\end{array}\right],
\end{equation*}
$A^{22}S^{111}=S^{221}A^{22}$ becomes $a^{12}=-a^{21}$ and
$a^{11}=a^{22}$. By~\rf{40} and the duality~\rf{10},
\begin{equation*}
N_{p_e}^{13}=\left[\begin{array}{cc}0&0\\0&
\left[\begin{array}{cc}0&\kappa\\-\kappa&0\end{array}\right]
\end{array}\right],
\end{equation*}
so that
\begin{equation*}
(N_{p_e}^{13})^*A^*+AN_{p_e}^{13}=
\left[\begin{array}{cc}0&0\\0&
-\kappa\left[
 \begin{array}{cc}2a^{12}&0\\0&2a^{12}
 \end{array}\right]
\end{array}\right],
\end{equation*}
so by~\rf{13}, the Hamiltonian generated by the new choice of $\Wred$
has the same form as~\rf{100} but with $I_1$ replaced by $I_1+2\kappa
a^{12}$. Thus, {\em the entire effect of the freedom to
choose~$\Wred$ is exactly that the constant $I_1$ may be arbitrarily
manipulated, as long as $\kappa\ne0$. The effective action of the
gauge group as it acts on the normal form is isomorphic to the
additive group of real numbers.}

\subsection{The normal form}

The splitting of $T_{p_e}P$ and subsequent analysis of the
linearization $dX_\Hxie(p_e)$ amounts to the construction of a normal
form for linearizations of relative equilibria. When equivariantly
moved through $G_\mue p_e$, this normal form can be combined with the
equivariant isotropic embedding
theorem~(\ct{MarsdenJE-1981.1},\ct{WeinsteinA-1997.1}), to provide a
$G_\mue$-equivariant normal form for the entire Hamiltonian system
near $G_\mue p_e$. I now construct this normal form, closely
following~\ct{PatrickGW-1995.1}.

Let $\tilde P$ be the symplectic manifold
\begin{equation*}\tilde P\equiv \Wred^0\times\Wred^1\times T^*G_\mue
\times G\mue,\end{equation*}
and define
\begin{equation*}\tilde p_e\equiv 0\oplus0\oplus0_\Id\oplus\mue.
\end{equation*}
Let the group $G_\mue$ act on $\tilde P$ by left translation on the
third factor. At the zero section of $T^*G_\mue$, there is the
splitting of the tangent space of $T^*G_\mue$ into vertical and
horizontal parts, and hence
\begin{equation}\lb{26}
T_{\tilde p_e}\tilde P=\Wred^0\oplus\Wred^1\oplus\mathfrak g_\mue\oplus
\mathfrak g_\mue^*\oplus T_\mue(G\mue).\end{equation} Together,~\rf{27}
and~\rf{26} yield a linear symplectomorphism between $T_{p_e}P$ and
$T_{\tilde p_e}\tilde P$, which extends by $G_\mue$-equivariance to a
vector bundle isomorphism, say $\Lambda^\prime:TP|(G_\mue
p_e)\rightarrow T\tilde P|(G_\mue \tilde p_e)$, which is symplectic on
each fiber. The equivariant isotropic embedding theorem then extends
$\Lambda^\prime$ to an equivariant symplectomorphism $\Lambda$ from a
neighborhood of $G_\mue p_e$ to a neighborhood of $G_\mue\tilde p$.
The relationship between $\Lambda^\prime$ and $\Lambda$ is
$\Lambda^\prime=T\Lambda$ on the domain of $\Lambda^\prime$.

The linearization $dX_\Hxie(p_e)$ is in a sense a first order
approximation of $X_\Hxie$ at $p_e$, and the quadratic Hamiltonian
associated to $dX_\Hxie(p_e)$ is a second order approximation to
$H_\xie$ at $p_e$. Together with the map $\Lambda$, this yields
a second order approximation to the Hamiltonian~$\tilde H$ in a
neighborhood of $G_\mue \tilde p_e$. Indeed, let $\tilde J^\xi$ be
the momentum mapping for the left action of right translations
$(g,h)\mapsto hg^{-1}$, so that
\begin{equation*}
\tilde J^\xi(\alpha_g)\equiv-\langle\alpha_g,X_\xi(g)\rangle,
\end{equation*}
where $X_\xi$ is the left invariant vector field generated by $\xi$.
Then 
\begin{equation}\lb{34}\tilde H=\Hnil-(\tilde J^{\xi_e}-\Hred^1)
+\Hred^0+\Hrsd+\Hrmd,
\end{equation} 
where
\begin{enumerate}
\item $\Hnil:\Wred^1\times T^* G_{\mu_e}\rightarrow\mathbb R$ is the 
extension of~\rf{33} to $T^* G_{\mu_e}$ by left invariance,
\item $\Hred^i:\Wred^i\rightarrow\mathbb R$ is 
$d^2 H_\mue(\bar p_e)|\Wred^i$, $i=1,2$,
\item $\Hrsd:G\mu_e\rightarrow\mathbb R$ by $\Hrsd(g\mu_e)\equiv-\langle g\mu_e,
\xi_e\rangle$, and
\item $d\Hrmd=0$ and $d^2\Hrmd=0$ on $G_{\mu_e} \tilde p_e$.
\end{enumerate}
Items (1)--(3) are definitions, while the content is in item~(4),
which follows since the two Hamiltonians
\begin{equation*}
\tilde H, \qquad
\Hnil-(\tilde J^{\xi_e}-\Hred^1)+\Hred^0+\Hrsd
\end{equation*}
both have relative equilibria at $\tilde p_e$ and have been
manipulated to have matched linearizations (and consequently matched
Hessians) at each point of $G_\mue\tilde p_e$. The point is that a
second order approximation to $\tilde H$ may be obtained by dropping
the remainder term $\Hrmd$ of~\rf{34}.

Now, consider $\tilde H$ as defined by~\rf{34}. The second term (the
one grouped within brackets) generates the one parameter group of
symplectomorphisms generated by $\xi_e$ through the
action~\rf{38}. Therefore, the first term $\Hnil$ and the second term
Poisson commute. Obviously, then, the first four terms of~\rf{34}
pairwise Poisson commute. The effect of the term $-\tilde J^\xie$ is
merely to generate the flow of the right invariant vector field
generated by $\xie$, and this corresponds to the ``fast'' evolution of
the relative equilibrium $\tilde p_e$, while the flows of the terms
$\Hred^0$ and $\Hrsd$ cannot cause evolution along $G_\mue\tilde
p_e$. So the slow evolution near $G_\mue\tilde p_e$ is to first order
dictated by the flow of $\Hnil$. 

As for the symmetries of $\tilde H$, since the Hamiltonian $H$ is $G$
invariant and $\Lambda^\prime$ is $G_\mue$-equivariant, the
Hamiltonian~\rf{34} is $G_\mue$-invariant, including the remainder
term. Also, when $\Hnil$ is extended to be a left invariant function
of $\Wred^1\times T^*G_\mue$, invariance under the normal form action
given by~\rf{11} becomes right invariance, i.e., for the $S^1$ part,
invariance under the action
\begin{equation}\lb{38}
a\cdot(x,\alpha_g)=\bigl(a^{-1}x,(TR_a)^*\alpha_g\bigr),
\end{equation}
where $R$ denotes right translation of $G_\mue$, and the functions
$\Hred^0$ and $\Hrsd$ are also invariant under~\rf{38}. However this
normal form symmetry does not in general extend to the remainder term
$\Hrmd$. The action~\rf{38} is counterclockwise or clockwise diagonal
rotation on the variables $(\pi_1,\pi_2)$ and $(x_1,x_2)$ in a way
that is a sense-preserving renormalization of the coadjoint action of
$\exp(\mathbb R\xi_e)$ on $(\pi_1,\pi_2)$. Since~\rf{38} is the action
of the Abelian $S^1$, this action can be reversed and thus assumed to
be counterclockwise; however, if this is done then the resulting
action may acquire the {\em opposite} sense to the coadjoint action of
$\exp(\mathbb R\xi_e)$ on $(\pi_1,\pi_2)$.

In concept the map $\Lambda$ is like a coordinate system near $G_\mue
p_e$; it is only that its range is a manifold rather than an open
subset of Euclidean space. Through its dependence on the splitting of
$T_{p_e}P$, $\Lambda$ depends on the choice of $\Wred$, so the effect
of varying $\Wred$ like varying a coordinate system near $G_\mue
p_e$. Thus, the freedom to vary $\Wred$ is like a gauge freedom. I
will call a particular choice of $\Wred$ \defemph{ a choice of gauge},
and say that, in view of Section~(3.1), the gauge group of $p_e$ is
$(\mathbb R,+)$. As shown in Section~(3.1), varying $\Wred$ itself is
a linear process, the effect of which on $\Hnil$ is easily found and
described. Yet, through the isotropic embedding theorem, $\Lambda$
depends on $\Wred$ in a nontrivial way.

This ability to cause very nontrivial ``coordinate changes'' while
easily calculating the concomitant effect on $\Hnil$ will enter
powerfully in the subsequent analysis. Particularly, one can imagine
choosing a gauge that ``simplifies'' $\Hnil$, making the analysis of
$\Hnil$ tractable and yielding information about the flow of the
original Hamiltonian system near $G_\mue p_e$. On the other hand, a
simple $\Hnil$ will in general be nongeneric, so that higher order
terms in $\Hrmd$ might destroy whatever delicate nongeneric structures
arise from $\Hnil$. In that case, one can focus on the identification
and perturbation of these nongeneric structures, which will occupy
smaller regions of phase space as ones attention is restricted more
nearly to $G_\mue \bar p_e$.

{\em Summary: The phase space
$\Pdrift$ on which some Hamiltonian $\Hdrift$ approximates to first order
the evolution of the system, when that system is perturbed from
$p_e$, is
\begin{gather*}
\Pdrift=\mathbb R^2\times T^*G_\mue,\\
\Hdrift(x,\alpha_g)=\langle g^{-1}\alpha_g,\xi_e\rangle+
\frac12N_{p_e}^{23}(\alpha_g,\alpha_g)+
\langle g^{-1}\alpha_g, N_{p_e}^{211}x\rangle.
\end{gather*}
The equations of motion for this system\footnote{\lb{52} For a
Hamiltonian system on the cotangent bundle of a Lie group $G$ of the
form $\alpha_g\mapsto
1/2h(\alpha_g,\alpha_g)+\langle\alpha_g,X\rangle$ where $h$ is a left
invariant metric on $G$ and $X$ is a left invariant vector field on
$G$, one has the equations of motion
\begin{equation*}
\Omega\equiv h^\sharp\pi+X(e),
\quad\pi\equiv g^{-1}\alpha_g,\quad g^{-1}\frac{dg}{dt}=\Omega,\quad
\frac{d\pi}{dt}=-\onm{coad}_\Omega\pi.
\end{equation*}}
are
\begin{equation*}
g^{-1}\frac{dg}{dt}=\Omega,\quad\frac{d\pi}{dt}=\onm{coad}_\Omega\pi,
\quad\frac{dx}{dt}=N^{131}_{p_e}\pi.
\end{equation*}
where
\begin{equation*}
\Omega\equiv\xi_e+(N^{23})^\flat\pi+N^{211}_{p_e}x.
\end{equation*}
There is a splitting $\mathfrak g_\mue=\{(\pi_1,\pi_2)\in \mathbb
R^2\} \oplus\{\mu^0\}$ such that $\Hdrift$ becomes,
after Poisson reduction by left translation of $G_\mue$,
\begin{equation}\lb{130}
\Hdrift=
\langle\mu^0,\xi_e\rangle+\frac12N_{p_e}^{230}(\mu^0,\mu^0)
+\frac12I_1({\pi_1}^2+{\pi_2}^2)
+\kappa(\pi_1x_1+\pi_2x_2).
\end{equation}
Here $N^{230}$ is a constant quadratic form, and $I_1$ and $\kappa$
are constants, and all of these are calculable from the nilpotent part
of the linearization of the relative equilibrium. The
Hamiltonian~\rf{130} is invariant under counterclockwise diagonal
rotation on the variables $(\pi_1,\pi_2)$ and $(x_1,x_2)$ and the
coadjoint action of $T_{\xi_e}$ on the variables $\mu^0$. The gauge
group of the approximation is $(\mathbb R,+)$, and the effect of the
gauge freedom is to exactly to undetermine $I_1$.}

\section{The mechanical case}

Here I consider the special case $G=SO(3)\times (S^1)^n$. Standard
identifications give $\mathfrak g=\mathbb R^3\times\mathbb
R^n=\mathfrak g^*$. In order that the momentum be nongeneric, assume
that the $SO(3)$ part of the momentum is zero. Reorient the system
(i.e. left translate it) so that $\mathfrak
g_\mue^1=\{(\pi_1,\pi_2,0)\}\subseteq so(3)$. Since the
action of $\onm{ad}_{\xi_e}$ fixes $\mathfrak g_\mue^1$, the $SO(3)$
part of $\xi_e$ is either parallel or antiparallel to to $\mathbf
k\in so(3)$, where $\mathbf k=(0,0,1)$. I drop the first term
of~\rf{130}, since the background motion of the relative equilibrium
$p_e$ itself will not be of interest, after which~\rf{130} becomes
\begin{equation*}\begin{split}
\Hdrift=\frac12I_1\bigl({\pi_1}^2+{\pi_2}^2\bigr)+\frac12I_2{\pi_3}^2
&+\kappa(\pi_1x_1+\pi_2x_2)\\
&\quad+\pi_3 \sum_{i=1}^na_ip_i+\sum_{i,j=1}^na_{ij}p_ip_j,
\end{split}\end{equation*}
where the variable labeling is
\begin{equation*}(\pi_1,\pi_2,\pi_3,p_1,\ldots,p_n)\in\mathfrak g^*=
\mathbb R^3\times\mathbb R^n,\end{equation*}
and $a_i$ and $a_{ij}$ are constants.

The whole of $(S^1)^n$ is ignorable; an Abelian reduction yields the 
phase space and Hamiltonian
(which by abuse of notation will have the same names)
\begin{gather}
\Pdrift=\mathbb R^2\times T^*SO(3)=\bigl\{\bigl((x_1,x_2)\in\mathbb R^2,
 A\in SO(3),(\pi_1,\pi_2,\pi_3)\in\mathbb R^3\bigr)\bigr\}
 \notag\\
\Hdrift=\frac12I_1({\pi_1}^2+{\pi_2}^2)+\frac12I_2{\pi_3}^2+
\kappa(\pi_1x_1+\pi_2x_2)+a\pi_3,\lb{41}
\end{gather}
where $a$ is constant. The normal form symmetry is counterclockwise
diagonal action of $S^1$ on the variables $(x_1,x_2), (\pi_1,\pi_2)$,
and the momentum mapping for the normal form symmetry is
\begin{equation*}J^{\mbox{\scriptsize nf}}=
-\frac12({x_1}^2+{x_2}^2)-\pi_3
\end{equation*}
An easy verification (see the footnote on page~\pageref{52}) gives
the equations of motion for~\rf{41} as
\begin{equation}
\frac d{dt}\left[\begin{array}{c}\pi_1\\\pi_2\\\pi_3\\x_1\\x_2
\end{array}\right]=
\left[\begin{array}{c}
I\pi_2\pi_3-\kappa x_2\pi_3+a\pi_2\\
-I\pi_1\pi_3+\kappa x_1\pi_3-a\pi_1\\
\kappa(\pi_1x_2-\pi_2x_1)\\
\kappa\pi_2\\-\kappa\pi_1\end{array}\right],\lb{54}\end{equation}
and
\begin{equation}\lb{55}
A^{-1}\frac{dA}{dt}=\left[
\begin{array}{c}I_1\pi_1+\kappa x_1\\I_1\pi_2+\kappa x_2\\I_2\pi_3+a\end{array}
\right]^\wedge,
\end{equation}
where $I\equiv I_2-I_1$ and as usual, for $v\in\mathbb R^3$,
\begin{equation*}v ^\wedge\equiv\left[\begin{array}{ccc}
\phantom{-}0\phantom{^3}&-v^3&\phantom{-}v^2\\
\phantom{-}v^3&\phantom{-}0\phantom{^2}&-v^1\\
-v^2&\phantom{-}v^1&\phantom{-}0\phantom{^2}
\end{array}\right].\end{equation*}

To obtain a first impression of~\rf{41} one can contrast it with the
nonresonant case. On the one hand, the phase space $\Pdrift$
of~\rf{41} has dimension~8; the symmetry group is $SO(3)\times S^1$.
Hence, the Marsden-Weinstein reduced systems are integrable, since
they have dimension $8-4-2=2$, and energy is a conserved quantity. We
shall see that homoclinic connections appear in the reduced spaces
of~\rf{41}, which are spheres in general, and that singular reduced
spaces (pinched spheres) also occur for some momentum values. On the
other hand, the {\em nonresonant\/} case (obtained from~\rf{41} by
setting $\kappa=0$ and dropping the variables $x_1$ and $x_2$) has the
same symmetry group but dimension~6 phase space; {\em it\/} is
integrable by its symmetry group alone, without the use of energy.
Dynamically, the nonresonant case is simpler, since every motion is a
relative equilibrium and the compact symmetry group rules out
homoclinic connections.

\subsection{Full reduction}
The Poisson manifold $\mathbb R^2\times\mathbb R^3=\{(x, \pi)\}$ may
be reduced by the $S^1$ normal form symmetry by use of the Hopf
variables (see, for example,~\ct{CushmanRHBatesLM-1997.1}, page~14; also
page~407 has a summary of references about singular reduction):
\begin{alignat*}{2}
w_1&\equiv 2(x_1\pi_2-x_2\pi_1),&\qquad w_2&\equiv 2(x_1\pi_1+x_2\pi_2),\\
w_3&\equiv({x_1}^2+{x_2}^2)-({\pi_1}^2+{\pi_2}^2),&\qquad
w_4&\equiv({x_1}^2+{x_2}^2)+({\pi_1}^2+{\pi_2}^2).
\end{alignat*}
This map is a quotient map for the $S^1$ normal form symmetry and it
has the semialgebraic image given by the subset of $\mathbb
R^5=\{(w_1,w_2,w_3,w_4,\pi_3)\}$ satisfying
\begin{equation*}
{w_1}^2+{w_2}^2+{w_3}^2-{w_4}^2=0,\quad w_4\ge 0.
\end{equation*}
The symplectic reduced spaces are the level sets of the two
Casimirs~$j_1$ and $j_2$ given by
\begin{gather}
\lb{70}
j_1\equiv({\pi_1}^2+{\pi_2}^2+{\pi_3}^2)^{\frac12}=\left({\pi_3}^2
+\frac12(w_4-w_3)\right)^\frac12,\\
\lb{71}j_2\equiv J^{\mbox{\scriptsize nf}}=
-\frac12\bigl({x_1}^2+{x_2}^2\bigr)-\pi_3=-\frac14(w_3+w_4)-\pi_3.
\end{gather}
Given $j_1$ and $j_2$,~\rf{70} and~\rf{71} can be used to eliminate 
$w_3$ and $w_4$, thereby obtaining the symplectic reduced spaces as the 
subsets of 
$\mathbb R^3=\{(w_1,w_2,\pi_3)\}$ given by
\begin{gather}
{w_1}^2+{w_2}^2=8({\pi_3}^2-{j_1}^2)(\pi_3+j_2)\lb{73}\\
w_4=({j_1}^2-2j_2)-({\pi_3}^2+2\pi_3)\ge0.\lb{7300}
\end{gather}

Various symplectic reduced spaces are obtained
by fixing various values of $j_1$ and $j_2$ , but 
the values of $j_1$ and $j_2$ are not arbitrary, since $j_1\ge0$ and
\begin{equation*}
j_1=({\pi_1}^2+{\pi_2}^2+{\pi_3}^2)^{\frac12}\ge-\pi_3=j_2+\frac12({x_1}^2+{x_2}^2)
\ge j_2.
\end{equation*}
Moreover, the singularities in the symplectic reduced spaces occur
where $x_1=x_2=\pi_1=\pi_2=0$, and here $j_1=|\pi_3|=\pm
j_2$. Equation~\rf{73} has solutions only over the intervals where the
cubic on the right hand side is nonnegative, and moreover, by~\rf{71},
$\pi_3\ge-j_2$. Putting all this together gives the bifurcation
diagram shown in~Figure~\rf{703}. As shown there, the
symplectic reduced spaces are the surfaces of revolution of
$\bigl(8({\pi_3}^2-{j_1}^2)(\pi_3+j_2)\bigr)^{\frac12}$ over the {\em
finite\/} interval where the cubic is positive. Thus, these spaces are
points, spheres, or, in the case that $j_2=-j_1$, they are topological
spheres with one conical singularity.
\begin{figure}[p]
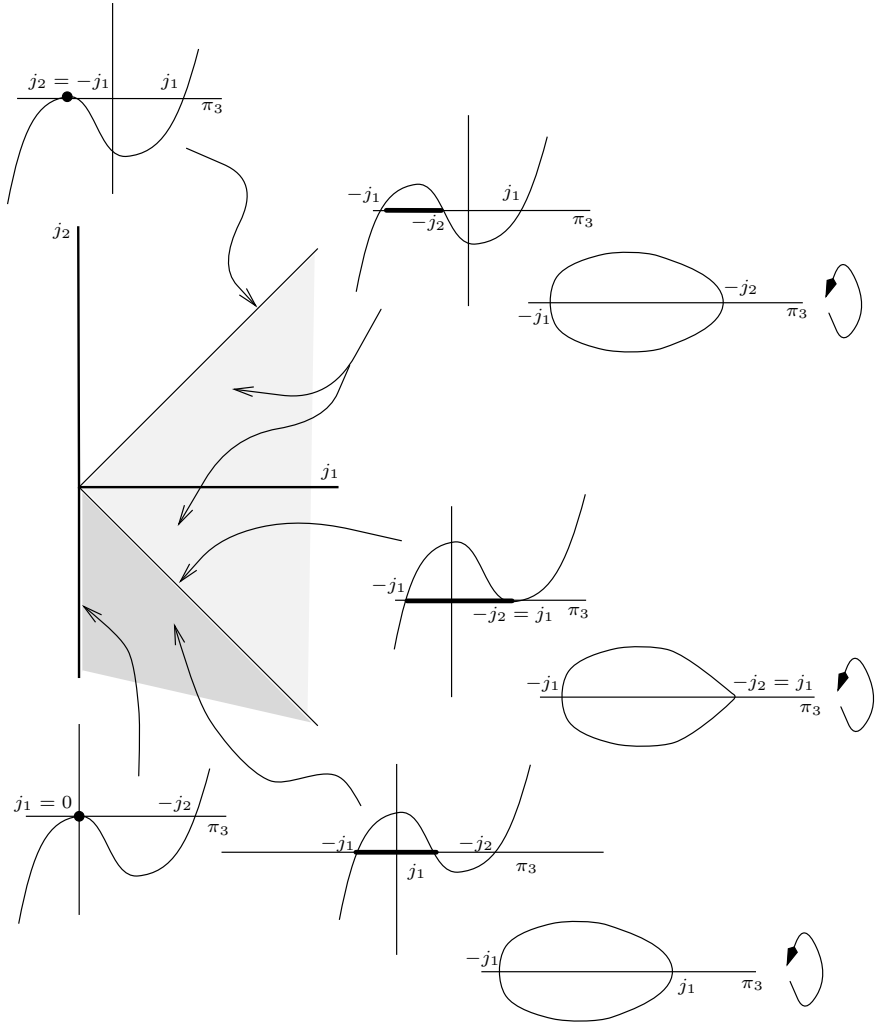
\lb{200}
\input fig1.tex\vspace*{.25in}
\caption{\label{703}\it\protect\footnotesize The bifurcation diagram
for the symplectic reduced spaces. Clockwise from the top, the reduced
spaces are points, spheres, pinched spheres, again spheres and finally
points. The finite intervals where the cubic is positive is the
thicker black line on the $\pi_3$ axis. One verifies that
Inequality~\rf{7300} is respected by showing that the concave down
quadratic $w_4=-({\pi_3}^2+2\pi_3)-{j_1}^2+2j_2$ is positive on those
intervals by checking the it is positive at the various endpoints of
the intervals.}
\end{figure}

In passing, I note that the symplectic volumes of the $-j_2=j_1>0$
reduced phase spaces are infinite, and of course the infinity is
concentrated at the singularity. In fact, by regarding two functions
$f(w_1,w_2)$ and $g(w_1,w_2)$ as functions on one of those reduced
phase spaces, lifting those functions to the Poisson phase space
$\{(x,\pi)\}=\mathbb R^2\times so(3)$, and then calculating
the Poisson bracket, one verifies that
\begin{equation*}
\{f,g\}=-4({\pi_3}^2+2j_2\pi_3+{j_1}^2)\left(
\frac{\partial f}{\partial w_2}\frac{\partial g}{\partial w_1}-
\frac{\partial f}{\partial w_1}\frac{\partial g}{\partial w_2}\right).
\end{equation*}
Consequently, the symplectic form on the reduced spaces in the coordinates
$(w_1,w_2)$ is
\begin{equation*}
\omega_{j_1,j_2}=\frac1{4({\pi_3}^2+2j_2\pi_3+{j_1}^2)}\,dw_1\wedge dw_2.
\end{equation*}
Switching to polar coordinates $w_1=r\cos\theta$ and $w_2=r\sin\theta$, the
symplectic volume of a small circle of radius $\delta$ about the origin
is
\begin{equation*}
A_\delta\equiv2\pi\int_0^\delta\frac{r}{4({\pi_3}^2+2j_2\pi_3+{j_1}^2)}\,dr
=\frac\pi2\int_0^\delta\frac{r}{(\pi_3-j_1)^2}\,dr
\end{equation*}
and upon substituting
\begin{equation*}
r^2=8({\pi_3}^2-{j_1}^2)(\pi_3-j_1)=8(\pi_3+j_1)(\pi_3-j_1)^2
\end{equation*}
that symplectic volume becomes
\begin{equation*}
A_\delta=4\pi\int_0^a\frac{\pi_3+j_1}{r}\,dr.
\end{equation*}
Since the numerator tends to $2j_1$ as $r$ becomes small, $A_\delta$ is
infinite.

The Hamiltonian~\rf{41} on the reduced spaces is, with the aid of
the Hopf variables and~\rf{70},
\begin{equation*}\begin{split}
\Hdrift&=\frac14I_1(w_4-w_3)+\frac12I_2{\pi_3}^2+
\frac\kappa2 w_2+a\pi_3\\
&=\frac12 I{\pi_3}^2+\frac\kappa2w_2+a\pi_3+\frac12I_1{j_1}^2
\end{split}\end{equation*}
Obviously this expression is simplified if $I_1=I_2$, since then
$I=I_2-I_1=0$. 
This choice of $I_1$ gives the first part of~\rf{41}
the same form as the kinetic energy of a spherical ball moving in
3-space, so I call this choice the \defemph{ spherical gauge}.

In the spherical gauge the flow lines of the reduced systems are
trivial to determine: they are the intersection of the planes
$(\kappa/2)w_2+a\pi_3=\mbox{constant}$ with the surfaces of revolution
in Figure~\rf{200}. In the case of a nonsingular reduced space of
nonzero dimension the flow is that is that of the flow on a 2-sphere with two
stable equilibrium points. For the singular reduced spaces
$-j_2=j_1>0$ the type of flow depends on whether or not the plane cuts
through the cone at the $\pi_3=j_1$ singularity
of~${w_2}^2=8(\pi_3+j_1)(\pi_3-j_1)^2$. That cone is
\begin{equation*}
|w_2|=4\sqrt{j_1}(j_1-\pi_3)\quad\Leftrightarrow\quad \pi_3=
j_1-\frac1{4\sqrt{j_1}}|w_2|,
\end{equation*}
and so there is an unstable equilibrium at the singularity if the line
$(\kappa/2)w_2+a\pi_3=aj_1$ passes through that cone, which is when
$a^2<4\kappa^2j_1$, and a stable equilibrium otherwise. The former
case is the flow on a 2-sphere with two stable equilibrium points and a
single homoclinic connection to a singular point, while in the latter
case this flow has two stable equilibrium points, one of which resides
at a singularity.

{\em For the remainder of this section I will assume the spherical gauge.}

The equilibria of~\rf{54} correspond to relative equilibria for just
the action of~$SO(3)$, since the space $\{(\pi,x)\}$ is the Poisson
phase space of $\Pdrift$ reduced by the $SO(3)$ action. These have a
greater significance than relative equilibria which also use the $S^1$
normal form symmetry, since the full unapproximated Hamiltonian has
the $SO(3)$ symmetry but it is only the truncated approximation which
has the $SO(3)$ symmetry \emph{and} the normal form symmetry. As
well, it is expedient to separate the equilibria of~\rf{54} at the
outset of the analysis. Setting the right side of~\rf{54} to zero
immediately gives $\pi_1=\pi_2=0$ and then $x_1\pi_3=x_2\pi_3=0$, so
there are the following solutions:
\begin{align}
&x=0,\quad\pi=0\lb{57}\\
&x=0,\quad\pi_1=\pi_2=0,\quad\pi_3\ne 0\lb{58}\\
&x\ne 0,\quad\pi=0\lb{59}
\end{align}
The set of relative equilibria given by~\rf{59} contains relative
equilibria equivalent under the $S^1$ normal form symmetry. The
$SO(3)$ generator corresponding to~\rf{59} is, by substitution
into~\rf{55}, $(\kappa x_1,\kappa x_2,0)$. The zero section of
$T^*SO(3)$ (i.e. the set $\pi=0$) corresponds to perturbations of the
original relative equilibrium $p_e$ such that the {\em perturbation}
has zero total angular momentum, and in the nonresonant case such
perturbations would imply no drift and no drift has been numerically
observed for such. For the resonant case, however, one expects regular
rotation around the axis $(\kappa x_1,\kappa x_2,0)$, which is
perpendicular to the generator $\xi_e$. Nonzero $x$ corresponds to
some ``reduced excitation''. Consequently, {\em at zero total angular
momentum the system drifts so that the generator of the relative
equilibrium moves along a fixed great circle at a rate dictated by a
reduced excitation.}

The relative equilibria~\rf{59} are the relative equilibria with
nongeneric momenta {\em for the system}~\rf{41}. Consequently, the
methods of~\ct{PatrickGW-1995.1} might be applied to them: the
drifting relative equilibria could themselves drift. I will not,
however, pursue this aspect here.

To obtain a list of nonequivalent relative equilibria (i.e. a list of
relative equilibria, no two of which are in the same $SO(3)\times S^1$
orbit), assume $A=\Id$ and discard Equation~\rf{55},
and, using the $S^1$ normal form symmetry, assume $\pi_2=0$ and
$\pi_1\ge0$, and when $\pi_1=0$ assume that $x_2=0$. Equating
the first of~\rf{54} with the $S^1$ infinitesimal generator of $s_e\in
\mathbb R$ gives
\begin{equation}\lb{62}
\left[
\begin{array}{c}
-\kappa x_2\pi_3+a\pi_2\\
\kappa x_1\pi_3-a\pi_1\\
\kappa(\pi_1x_2-\pi_2x_1)\\
\kappa\pi_2\\-\kappa\pi_1\end{array}
\right]
=\left[\begin{array}{c}-s_e\pi_2\\s_e\pi_1\\0\\-s_ex_2\\s_ex_1
\end{array}\right].\end{equation}
Putting $\pi_2=0$ into the fourth component gives $s_ex_2=0$. Now
$s_e=0$ corresponds to the $SO(3)$ relative equilibria~\rf{57}--\rf{59},
and in~\rf{57} and~\rf{58}, we have $x_2=0$, while in~\rf{59} we can
assume $x_2=0$ and $x_1>0$ by the normal form symmetry. So if $s_e=0$
then we assume $x_2=0$, while if $s_e\ne0$ then
$x_2=0$ anyway by $s_ex_2=0$. Putting $x_2=\pi_2=0$ into~\rf{62} then
gives
\begin{equation*}
-\kappa\pi_1=s_ex_1,\qquad 
\kappa x_1\pi_3-a\pi_1=s_e\pi_1,
\end{equation*}
and these equations are easily solved to obtain the following list of
nonequivalent relative equilibria:
\begin{align}
&\pi_1=\pi_2=x_1=x_2=0,\eta_e=(a+s_e)\mathbf k,\;\pi_3\in\mathbb R;\lb{64}\\
\begin{split}&\pi_2=x_2=0,
\pi_3=\frac{\pi_1(ax_1-\kappa\pi_1)}{\kappa {x_1}^2},
s_e=-\frac{\kappa\pi_1}{x_1},\\
&\qquad\qquad\qquad\eta_e=\kappa x_1\mathbf i+(a+s_e)\mathbf k,\;\pi_1\in
\mathbb R,x_1>0.\end{split}\lb{65}
\end{align}
Here the $SO(3)$ generators $\eta_e$ of these relative equilibria have
been determined by comparison of~\rf{55} and the infinitesimal
generator of the $SO(3)$ action of $\Pdrift$ at the relative
equilibrium. In~\rf{64} $\eta_e$ and $s_e$ are not unique due to the
presence of isotropy in the phase space $\Pdrift$. For later use, the
characteristic polynomials of the linearizations of these relative
equilibria are, respectively,
\begin{equation*}
p_1\equiv x^2(x^2+|\eta_e|^2)
\bigl(x^4+(a^2+2as_e+2s_e^2-2\pi_3\kappa^2)x^2+
(s_e^2+as_e+\pi_3\kappa^2)^2\bigr)
\end{equation*}
and
\begin{equation*}
p_2\equiv
\frac1{{x_1}^2}x^4(x^2+|\eta_e|^2)\bigl({x_1}^2x^2+\kappa^2{x_1}^4
+(ax_1-2\kappa\pi_1)^2\bigr).
\end{equation*}

For an understanding of the relative equilibria and the flow in
general it is important to understand which relative equilibria reside
on which reduced spaces. For~\rf{64}, $j_1=|\pi_3|$ and $j_2=\pi_3$,
so for $\pi_3\le0$ one has $j_1=j_2$, which corresponds to pointlike
reduced spaces, while for $\pi_3>0$ one has $j_1=-j_2$, which
corresponds to the singular reduced spaces, and the relative
equilibria occupy the singularities. As for~\rf{65},
substituting~$\rf{65}$ into $j_1$ and $j_2$ and eliminating $\pi_1$
gives
\begin{multline}\lb{201}
(3{x_1}^4+8j_2{x_1}^2+4{j_2}^2-4{j_1}^2)^2\kappa^2\\
+4a^2{x_1}^2({x_1}^2+2j_2-2j_1)({x_1}^2+2j_2+2j_1)=0
\end{multline}
so it is a matter of solving this quartic in ${x_1}^2$ for its
positive roots. The pointlike reduced spaces corresponding to $j_1=0$
are occupied by the relative equilibria of~\rf{65} after putting
$\pi_1=0$, corresponding to the solutions $j_1=0$, $x_1=\sqrt{-2j_2}$
of~\rf{201}.

\subsection{Partial reduction by the normal form symmetry}
In this section I symplectically reduce the Hamiltonian system by the
normal form symmetry at a value, say $\sigma$ of the momentum mapping
$J^{\mbox{\scriptsize nf}}$. My objective is an interpretation of the
reduced system via the equations of motion on the reduced space.

Use the notation $r^2={x_1}^2+{x_1}^2$. 
The $\sigma$-level of the momentum map $J^{\mbox{\scriptsize nf}}$ is
\begin{equation*}
\pi_3=-\sigma-\frac12r^2.
\end{equation*}
Below, when $x$ appears in the context of a vector in $\mathbb R^3$,
it is as $(x_1,x_2,0)$.
The map
\begin{equation}\lb{80}
\Psi(A,\pi,x)\equiv(y,\alpha,z)\equiv
\bigl(A\mathbf k,A(\pi-(\pi\cdot\mathbf k)\mathbf k),Ax\bigr)
\end{equation}
is a projection of the phase space $T^*SO(3)\times\mathbb R^2$ to the
Whitney direct sum $T^*S^2\oplus TS^2$. The restriction of this map to
the $\sigma$ level set of $J^{\mbox{\scriptsize nf}}$ is clearly a
quotient map for the normal form action, and thus the symplectic
reduced space is $T^*S^2\oplus TS^2$. The action of $SO(3)$ on 
$T^*S^2\oplus TS^2$ becomes
\begin{equation*}
A\cdot(y,\alpha,z)=(Ay,A\alpha,Az)
\end{equation*}

To find the reduced vector field at $(y,\alpha,z)\in T^*S^2\oplus TS^2$,
first evaluate the original Hamiltonian vector field, namely
\begin{alignat}{2}
&
\Omega\equiv I_1\pi+I((\pi\cdot\mathbf k)+a)\mathbf k,&\qquad
&A^{-1}\frac{dg}{dt}=\Omega^\wedge,\nonumber\\
&\frac{d\pi}{dt}=\dot\pi=\pi\times\Omega+a\mathbf k,
&\qquad\qquad&\frac{dx}{dt}=\dot x=-\kappa\mathbf k\times\pi\lb{81}
\end{alignat}
at the point
\begin{equation}
\pi= A^{-1}\alpha-\left(\sigma+\frac12r^2\right)\mathbf k,\quad
x=A^{-1}z,\lb{82}
\end{equation}
where $A$ is chosen so that $A\mathbf k=y$, and then apply the
derivative of~\rf{80}, which is
\begin{equation}\begin{split}
\lefteqn{T\Psi(A,\pi,x,\Omega,\dot\pi,\dot x)}\\
&\equiv
\dbyd t0\Psi\bigl(A\exp(\Omega^\wedge t),\pi+t\dot\pi,x+t\dot x\bigr)\\
&=A\cdot\bigl(\Omega\times\mathbf k,
 \Omega\times(\pi-(\pi\cdot\mathbf k)\mathbf k)+\dot\pi
 -(\dot\pi\cdot\mathbf k)\mathbf k,
 \dot x+\Omega\times x\bigr).
\end{split}\lb{83}\end{equation}
Thus it is a matter of substituting~\rf{81} and~\rf{82} into~\rf{83}.
Without care the calculation can be onerous; however, it is
straightforward, and yields the equations of motion
\begin{gather}
\frac{dy}{dt}=(I_1\alpha+\kappa z)\times y,\lb{841}\\
\frac{d\alpha}{dt}=-\pi_3(I_1\alpha+\kappa z)\times y
 +\kappa z\times \alpha,\lb{842}\\
\frac{dz}{dt}=\big(\kappa\alpha-(I_2\pi_3+a)z\bigr)\times y
 +I_1\alpha\times z,\lb{843}
\end{gather}
where $\pi_3$ stands for $-(\sigma+|z|^2/2)$. To these equations
must be added the constraints $|y|=1$, $\alpha\cdot y=0$ and $z\cdot
y=0$.

I want to impose the viewpoint that there is a ``particle'' at~$y$
having direction~$z$. To do this, I use~\rf{841} to replace $\alpha$
with $dy/dt$ in~\rf{842} and~\rf{843} while writing equation~\rf{842}
as a second order equation in~$y$. I also use the standard
Levi-Cevita connection of $S^2$ for the time derivatives. Thus,
Equation~\rf{842} becomes
\begin{equation*}\begin{split}
\frac{d^2y}{dt^2}&=(I_1\alpha+\kappa z)\times\frac{dy}{dt}+
\left(I_1\frac{d\alpha}{dt}+\kappa\frac{dz}{dt}\right)\times y\\
&=\bigl(I_1\alpha+\kappa z\bigr)\times\bigl((I_1\alpha+\kappa z)
 \times y\bigr)\\
&\qquad\mbox{}+I_1\pi_3(I_1\alpha+\kappa z)-\kappa\bigl(\kappa\alpha
 -(I_2\pi_3+a)z\bigr)\\
&=-\left|\frac{dy}{dt}\right|^2y+I_1\pi_3y\times\frac{dy}{dt}
-\kappa\left(\frac\kappa{I_1}
 \left(y\times\frac{dy}{dt}-\kappa z\right)-(I_2\pi_3+a)z\right),
\end{split}\end{equation*}
and similarly with Equation~\rf{843}, so that
\begin{gather}
\frac{\nabla^2y}{dt^2}=\frac{{I_1}^2\pi_3-\kappa^2}{I_1}y\times\frac{dy}{dt}
+\kappa\left(\frac{\kappa^2}{I_1}+I_2\pi_3+a\right)z\lb{85}\\
\frac{\nabla z}{dt}=\frac\kappa{I_1}\frac{dy}{dt}
 +\left(\frac{\kappa^2}{I_1}+I_2\pi_3+a\right)y\times z.\lb{86}
\end{gather}

Some aspects of the particle become apparent by replacing 
$\alpha$ with $dy/dt$ in the total energy and angular momentum.
Using~\rf{82}, the total energy is
\begin{equation}\begin{split}
H&=\frac12|\alpha|^2+\frac12{\pi_3}^2+a\pi_3+\kappa z\cdot\alpha\\
&=\frac1{2I_1}|I_1\alpha+\kappa z|+\frac12{\pi_3}^2+a\pi_3
-\frac{\kappa^2}{2I_1}|z|^2\\
&=\frac{1}{2I_1}\left|\frac{dy}{dt}\right|^2
 +\frac{1}{2I_2}\left(\frac{\kappa^2}{I_1}+I_2\pi_3+a\right)^2
 +\frac{\kappa^2\sigma}{I_1},\end{split}\lb{93}
\end{equation}
so the internal energy of the particle is
\begin{equation*}
E_{\mbox{\scriptsize int}}
=\frac{1}{2I_2}\left(\frac{\kappa^2}{I_1}+I_2\pi_3+a\right)^2.
\end{equation*}
The total angular momentum, or the $SO(3)$-momentum-mapping is,
using~\rf{82},
\begin{equation}
J=A\pi=\alpha-\pi_3\mathbf k=\frac1{I_1}y\times\frac{dy}{dt}
-\frac{\kappa}{I_1}z-\pi_3y,\lb{92}
\end{equation}
and below I replace $z$ with the tangent to the sphere part of the angular 
momentum, which by~\rf{92} is 
\begin{equation*}
L\equiv-\frac\kappa{I_1}z,
\end{equation*}
while the total angular momentum attributed to the particle must
be $L+\pi_3y$.
By the first terms of both~\rf{93} and~\rf{92}, the particle should be 
viewed as having mass $1/I_1$.

From~\ct{JacksonJD-1975.1} the classical nonrelativistic equations of 
motion on $|y|=1$
of a particle of mass $m$,
charge $Q$, and gyromagnetic ratio $\Gamma$, and magnetic 
moment $\mathbf m=\Gamma L$, are
\begin{gather}
m\frac{\nabla^2y}{dt^2}=
\frac QcB\times \mathbf v+\nabla(\mathbf m\cdot B)-(\nabla\cdot B)\mathbf m,\lb{87}\\
\frac{\nabla L}{dt}={\mathcal I}+\Gamma L\times B,\lb{88}
\end{gather}
where $c$ is the velocity of light and $\mathbf v$ is the particle's velocity.
Here $\mathcal I$ is compensates for the inductive forces that would, 
were $\mathbf m$ to be constant, be required to maintain the constant current 
that generates~$\mathbf m$ itself.
One can fit equations~\rf{85} and~\rf{86} to these equations as follows.
If one takes $B=y$, comparison of the first terms of~\rf{86} and~\rf{88}
gives the inductive term 
\begin{equation}
{\mathcal I}=-\frac{\kappa^2}{{I_1}^2}\frac{dy}{dt},\lb{95}
\end{equation}
while the second term of~\rf{86} matches the second term of~\rf{88} 
if the gyromagnetic ratio is
\begin{equation}
\Gamma=\frac{{I_1}^2}{\kappa^2}\left(\frac{\kappa^2}{I_1}+I_2\pi_3
+a\right).\lb{89}
\end{equation}
Matching the second term of~\rf{85} with the third term of~\rf{87}
implies $\Gamma$ is $\kappa^2/3I_1$ times the left side of~\rf{89}. So
choose a gauge so that $I_1=\kappa^2/3$. 
Then matching the first term of~\rf{85} with the first term of~\rf{87}
yields
\begin{equation*}
\frac Qc=\pi_3-\frac{\kappa^2}{{I_1}^2}.
\end{equation*}
Finally, I note that the inductive term~\rf{95} is exactly
what is required to balance the rate of change of internal energy with the
work done on the particle through its magnetic moment via the last term
of~\rf{85}.

\section{Example: Two coupled rods}

In this section I illustrate and verify the above theory using
numerical simulations of the system of two axially symmetric rods
which are joined by a frictionless ball-and-socket joint. In the
center of mass frame this system can be cast as a geodesic flow with
configuration space $SO(3)^2=\{(A_1,A_2)\}$ and kinetic
energy metric
\begin{equation*}
L=\frac{1}{2}\left[\begin{array}{cc}{\Omega_1}^t&
 {\Omega_2}^t\end{array}\right]J(A)\left[\begin{array}{c}\Omega_1\\ \Omega_2
 \end{array}\right],
\end{equation*}
where $A=[A^{ij}]={A_1}^tA_2$, $\Omega_1$ and $\Omega_2$ are body 
referenced (left translation) angular velocities, and
\begin{equation*}J(A)=\left[\begin{array}{cccccc}
1&0&0&-\beta A^{22}&\phantom{-}\beta A^{21}&0\\
0&1&0&\phantom{-}\beta A^{12}&-\beta A^{11}&0\\
0&0&\alpha&0&0&0\\
-\beta A^{22}&\phantom{-}\beta A^{12}&0&1&0&0\\
\phantom{-}\beta A^{21}&-\beta A^{11}&0&0&1&0\\
0&0&0&0&0&\alpha\end{array}\right].
\end{equation*}
Here $0<\alpha<2-2\beta$ is a parameter which increases with the
diameter of the rods and $0\le\beta<1$ is a parameter which measures
the degree of coupling between the rods; $\beta$ is zero if the joint
lies at the mutual centers of mass of the rods. The continuous
symmetry group of this system is $SO(3)\times (S^1)^2$, which acts on
the configuration space $SO(3)^2$ by
\begin{equation*}
(B,\theta _1,\theta _2 )\cdot(A_1,A_2)=\bigl(BA_1\exp(-\theta_1
\mathbf k^\wedge),BA_2\exp(-\theta_2\mathbf k^\wedge)\bigr).
\end{equation*}

All relative equilibria for this system are explicitly known, as are
their formal stability; a complete list can be found
in~\ct{PatrickGW-1991.1}. The relative equilibria of interest here are
the phase space points $\bigl(\Id,\exp(\theta\mathbf
j^\wedge),\Omega_1,\Omega_2\bigr)$ parameterized
by~$t_1,t_2,\theta\in\mathbb R$ such that $t_1\ne0,t_2\ne0,
0<\theta<\pi$, where
\begin{equation*}
\Omega_1=t_1\mathbf i-\kappa^{t_1t_2\theta}_1\mathbf k,
 \Omega_2=t_2\mathbf i+\kappa^{t_2t_1\theta}_1\mathbf k
\end{equation*}
and where
\begin{equation*}
\kappa^{t_1t_2\theta}_\gamma\equiv\frac{(t_1\cos \theta-t_2)(\beta t_2-\gamma t_1)}
 {\alpha t_1\sin \theta}.
\end{equation*}
The corresponding generators
$(\Omega,\sigma_1,\sigma_2)\in so(3)\times\mathbb R^2$ are given by
\begin{equation*}
\sigma_1=\kappa^{t_1t_2\theta}_{1-\alpha},\sigma_2=-
 \kappa^{t_2t_1\theta}_{1-\alpha},\Omega=t_1\mathbf i+\frac{t_1\cos \theta-t_2}
 {\sin \theta}\mathbf k.
\end{equation*}
By calculating the linearizations, one sees that a 1:1 group reduced
resonance with zero total-angular-momentum can be arranged by setting
\begin{equation*}
\beta=\frac{2t_1t_2}{{t_1}^2+{t_2}^2},
\quad 3{t_1}^4+{t_1}^2{t_2}^2(4\cos^2\theta-10)+3{t_2}^4=0,
\end{equation*}
and for the work presented in this section I have chosen to perturb the
particular relative equilibrium obtained by setting
\begin{equation}\lb{900}
\theta=\frac\pi3,t_1=1,t_2=\frac{1+\sqrt 5}{2},
\alpha=\frac12,\beta=\frac2{\sqrt 5}.
\end{equation}

\subsection{Calculations for the coupled rod system}
To compare the dynamics of the drift system to the dynamics of the coupled rod
system near the relative equilibrium~\rf{900}, it is necessary to
calculate the splitting~\rf{1} or~\rf{2} in this special case. Here is
one general way, inspired by the proof of Theorem~3.1.19
of~\ct{AbrahamRMarsdenJE-1978.1}:
\begin{enumerate}\renewcommand{\labelenumi}{(\alph{enumi})}
\item Calculate the semisimple part $S_{p_e}$ and the nilpotent part
$N_{p_e}$ of the linearization $dX_{H_\xie}(p_e)$.
\item Since $S_{p_e}$ is semisimple it has a basis of eigenvectors,
some complex and some real, and the complex eigenvectors may be
grouped into complex conjugate pairs. Taking the real and imaginary
parts of one eigenvector in each pair, and then including the real
eigenvectors yields a basis $\mathcal B$ of $T_{p_e}P$ in which
$S_{p_e}$ is {\em skew symmetric}.

\item Let the matrix of $\omega(p_e)$ with respect to the basis $\mathcal
B$ be $-W$. Since $S$ is infinitesimally symplectic and skew,
$0=S^tW+WS=-SW+WS$, so $S$ and $W$ commute. By going to a basis of
eigenvectors of the positive matrix $-W^2$, find a symmetric, positive
square-root $B$ of $-W^2$ such that $B$ and $S$ commute, and set
$J=WB^{-1}$. Using the basis $\mathcal B$, regard $J$ as its
corresponding operator on $T_{p_e}P$, and regard $B$ as the bilinear
form on $T_{p_e} P$ corresponding to $(x,y)\mapsto x^tBy$. Then
\begin{equation*}
\omega(p_e)(v,w)=B(Jv,w),\quad J^2=-\Id,
\end{equation*}
and $J$ commutes with $S_{p_e}$.

\item Set $Z=J(\mathfrak g_\mue)$, so $Z$ is a $S_{p_e}$ invariant
Lagrangian complement to $\mathfrak g_\mue$ in $\mathfrak g_\mue\oplus
Z$. Choose an $\onm{Ad}$-invariant complement $\mathfrak b$ to $\mathfrak
g_\mue$, and set $\Wred=(\mathfrak g_\mue p_e\oplus Z\oplus\mathfrak
bp_e)^{\omega\perp}$. Then $\Wred$ is symplectic, $S_{p_e}$ invariant
and is contained in $\onm{ker} dJ(p_e)$, since $\Wred\subseteq
\mathfrak gp_e^{\omega\perp}$. The subspace $\Wred$ complements
$\mathfrak g_\mue$ in $\onm{ker}dJ(p_e)$ since $\mathfrak
g_\mue\subseteq \Wred^{\omega\perp}$ and $\Wred$ is symplectic.
\end{enumerate}
In the case of two coupled rods, further refinement of the splitting
of $\Wred$ to resonant and nonresonant parts is not necessary, since
the reduced spaces for zero total-angular-momentum are two
dimensional; there is no ``nonresonant'' part of the phase space.
Final adjustments to
achieve the spherical gauge are easily arranged using~\rf{12} and~\rf{13}. The end result of all this is a
basis of the tangent space at the relative equilibrium which reflects
the splitting~\rf{2}, a basis of the Lie algebra $so(3)\times\mathbb
R^3$, reflecting the splitting~\rf{102} into resonant and nonresonant
parts, the information
\begin{equation}\lb{3000}
\kappa=.9115064,\quad I_1=I_2=4.321619,
\end{equation}
and the information that~\rf{900} is a  ``$+$''~type
relative equilibrium. 

Further to calculating the splitting~\rf{1} or~\rf{2}, means must be
found to translate initial conditions of the drift system to
perturbations (i.e. initial conditions near~\rf{900}) of the coupled
rod system. For this a map was used, say $\Psi$, of the tangent space
$T_{p_e}P$ into phase space $P$ such that $T_{p_e}\Psi$ is the
identity at $p_e$. The map $\Psi$ only injected initial conditions of
the drift system into perturbations of~\rf{900} to first order; actual
correspondence could not be easily achieved due to the implicit nature of
the coordinates provided by the isotropic embedding theorem.
Moreover, $\Psi$ did not match the momentum of the drift
system exactly to a momentum perturbation of the actual system, while
practice  indicated that an exact match of these momenta is
important. For the comparisons just below, given an initial condition
of the drift system with a particular momentum, an initial condition
of the coupled rod system with matching momentum was obtained as
$\Psi$ of another nearby initial condition of the drift system. That
nearby initial condition of the drift system was obtained by a slight
(second order) iterative refinement of the original initial condition
of the drift system.

By the conventions of Section~4, and since the relative
equilibrium~\rf{900} is of ``$+$''~type, the vector
$(\pi_1,\pi_2,\pi_3)=\mathbf k$ is parallel to the normalized rotation
vector $\xi_e$ of the coupled rod system. Also, the quantity $A\in
SO(3)$ of the drift system corresponds to the drifted orientation of
the coupled rod system. Thus, the quantity $A\mathbf k$ of the drift system
corresponds to the drifting rotation vector of the coupled rod
system. Now, the rotation vector of the relative equilibrium~\rf{900}
is a constant linear sum of the locations of the rods, exactly because
the rods when in the relative equilibrium rotate around that vector: if
one defines
\begin{equation*}
\tau\equiv \frac{1}{\sin\theta}\left(t_1A_1\mathbf k -t_2A_2\mathbf k\right)
\end{equation*}
then $\tau=\xi_e$ at the relative equilibrium~\rf{900} and also at any
reorientation of~\rf{900}. Thus, for a perturbation of~\rf{900},
$\tau$ approximates the drifting rotation vector and the prediction of
the preceding theory is that
\begin{equation*}
\bar\tau\equiv\frac\tau{|\tau|}\approx A\mathbf k.
\end{equation*}
The numerical  verifications I have  undertaken consist of predictions
of the drift system for the motion of $A\mathbf  k$ and the comparison
of these predictions with evolution  of the the unit vector $\bar\tau$
calculated  via  a symplectic integration of   the coupled rod system.
The particular algorithm  used was  a (implicit) Riemannian   leapfrog
algorithm~(\ct{PatrickGW-1996.1}).

\begin{figure}\setlength{\unitlength}{1in}\centering
\begin{picture}(4.5,2)
\put(1.,1){\makebox(0,0){\epsfbox{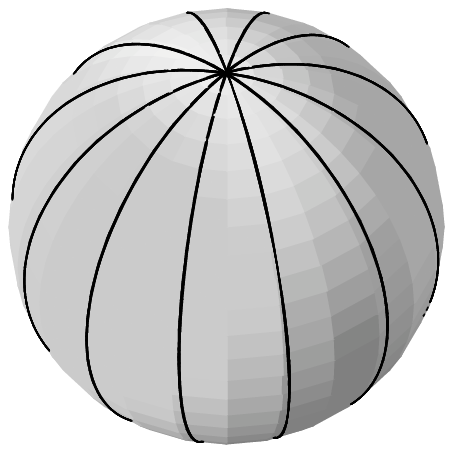}}}
\put(3.4,1){\makebox(0,0){\epsfbox{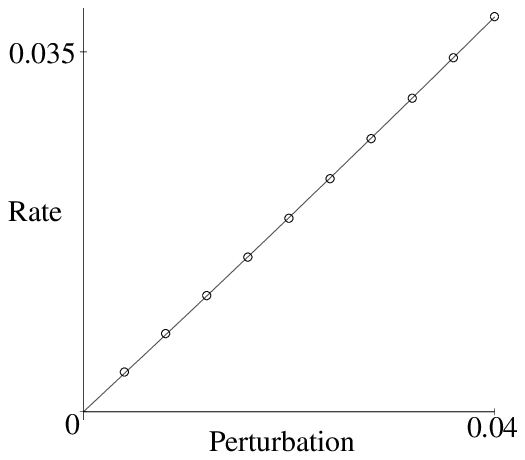}}}
\end{picture}
\caption{\lb{2500}\protect\footnotesize\it Perturbations within
zero-total-angular momentum. Left: the paths traced out by $\bar\tau$,
resulting from $6$~initial conditions within zero
total-angular-momentum. Right: the rotation rates of $\tau$ vs. the
magnitude of the perturbation $|x|$.}
\end{figure}
\subsection{Numerical results}
\subsubsection{First comparison: zero total-angular-momentum}
By the results in Section~(4.1), the part of the 
drift phase space corresponding
to perturbations with zero total-angular-momentum is occupied by
$SO(3)$ relative equilibria, and the motion of $A\mathbf k$ from initial
condition $\pi=0$ and $x$ arbitrary is that of uniform rotation at
angular frequency $\kappa|x|$ along a great circle through $\mathbf k$
and perpendicular to $(x_1,x_2,0)$. Translated to the coupled rod
system, the prediction is that $\bar\tau$ undergoes regular rotation
in a great circle through $\bar\tau(0)$, the particular great circle
regularly rotating as $x$ is rotated. The left of Figure~\rf{2500}
shows the motion of $\tau$ as $x$ is varied from $0$ to $0$ through
$5\pi/6$ radians by increments of $\pi/6$, so the coupled rod system
conforms to this prediction. As the magnitude of $x$ is varied the
great circles should be traced out at an angular rotation rate $\kappa
|x|$, so $\kappa$ can be determined by plotting that rotation rate
against $|x|$; this is done on the right of Figure~\rf{2500}. The
rotation rates fit well to the curve
\begin{equation*}
\mbox{Rate}= .9284|x|+.7750|x|^2,
\end{equation*}
so from the simulation $\kappa$ is $.9284$, which is $1.8\%$ off the
calculated value of $.9115064$ already displayed in~\rf{3000}.

\subsubsection{Second comparison: a stable relative equilibrium of the 
drift system}
For another comparison, I examined the coupled rod system for
perturbations of~\rf{900} corresponding to being near one of the
stable relative equilibria in the list~\rf{65}. For perturbations
of~\rf{900} corresponding to the values
\begin{equation}\lb{1001}
\pi_1=.001,\quad x_1=.04
\end{equation}
in~\rf{65} the predictions of the drift system are that $\tau$ moves
as small periodic oscillation of the full reduction of the drift
system superimposed on the periodic motion $\exp(\eta_e^\wedge t)\mathbf
k$, where $\eta_e$ is the $SO(3)$ generator in the list~\rf{65}. The
two predicted frequencies are the linearized frequency of the
corresponding equilibrium on the reduced space of the drift system
and the rotation
frequency of $\eta_e$, which are by substitution of~\rf{1001}
into~\rf{65} and/or~\rf{201}, respectively,
\begin{equation}\lb{2004}
.04477,\quad .05837.
\end{equation}
Again by substitution of~\rf{1001} into~\rf{65}, the motion of $\tau$
should be such that its projection onto the unit vector
$\eta_e/|\eta_e|$ is near $-.5300$.

\begin{figure}[p]\setlength{\unitlength}{1in}\centering
\begin{picture}(4.5,2)
\put(1.,1){\makebox(0,0){\epsfbox{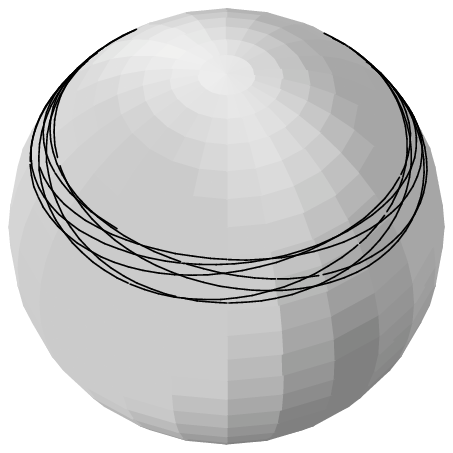}}}
\put(3.4,1){\makebox(0,0){\epsfbox{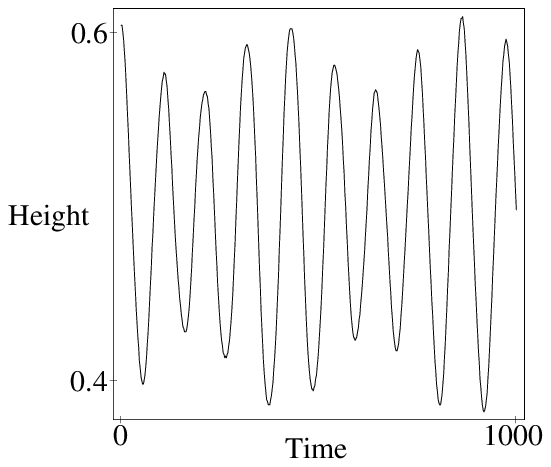}}}
\end{picture}

\vspace*{.05in}\begin{picture}(4.5,2)
\put(2.25,1){\makebox(0,0){\epsfbox{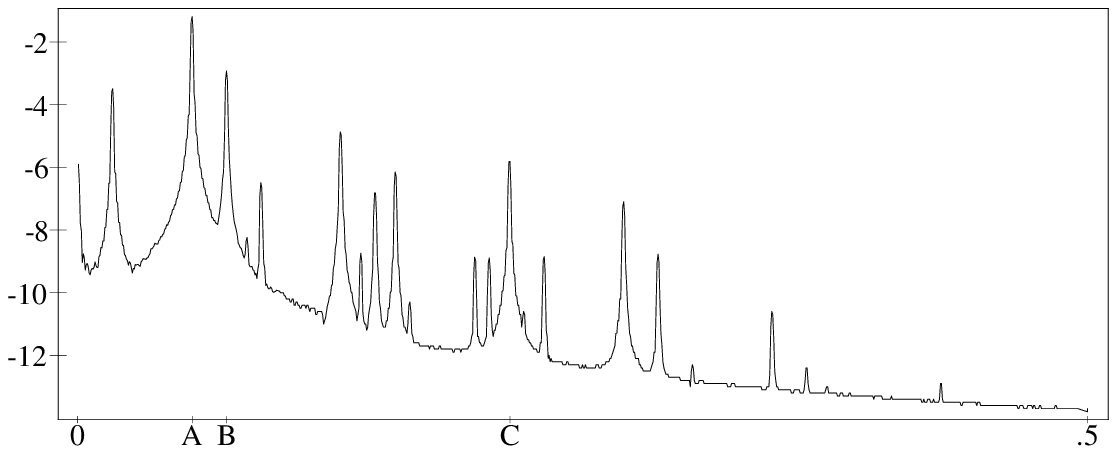}}}
\end{picture}

\vspace*{.15in}\footnotesize\begin{tabular}{rrrrrr|rrrrrr}
Power&Freq.&Freq.&$z_1$&$z_1$&$z_3$&
Power&Freq.&Freq.&$z_1$&$z_1$&$z_3$\\
\hline
1.195 & .04449 & .04449 & 1 & 0 & 0 &
8.766 & .2266 & .2266 & 0 & 1 & 1 \\
2.935 & .05791 & .05791 & 0 & 1 & 0 &
8.848 & .1822 & .1822 & -1 & 1 & 1 \\
3.485 & .01342 & .01342 & -1 & 1 & 0 &
8.874 & .1549 & .1553 & 1 & -1 & 1 \\
4.871 & .1024 & .1024 & 1 & 1 & 0 &
8.893 & .1607 & .1603 & 1 & 2 & 0 \\
5.812 & .1687 & .1687 & 0 & 0 & 1 &
10.3 & .1296 & .1292 & -1 & 3 & 0 \\
6.146 & .1239 & .1243 & -1 & 0 & 1 &
10.58 & .1741 & .1737 & 0 & 3 & 0 \\
6.485 & .07133 & .07133 & -1 & 2 & 0 &
10.63 & .2711 & .2711 & 1 & 1 & 1 \\
6.817 & .1162 & .1158 & 0 & 2 & 0 &
12.25 & .2401 & .2401 & -1 & 2 & 1 \\
7.102 & .2132 & .2132 & 1 & 0 & 1 &
12.36 & .2846 & .2846 & 0 & 2 & 1 \\
8.24 & .06596 & .06634 & -1 & -1 & 1 &
12.91 & .3371 & .3375 & 0 & 0 & 2 \\
8.742 & .1104 & .1108 & 0 & -1 & 1 &
12.98 & .2926 & .293 & -1 & 0 & 2 \\
\end{tabular}
\vspace*{.15in}\footnotesize
\caption{\lb{1000}\protect\footnotesize\it
Some results of the simulation of the
coupled rod system for a perturbation corresponding to a elliptic
relative equilibrium of the drift system. Top left: the evolution of
the rotation vector projected onto the unit sphere. Top right: the
time evolution of the height corresponding to the top left. Middle:
the power spectrum of a signal derived from the drifting motion; note
the vertical logarithmic scale. Bottom table: the peaks in the power
spectrum, sorted by power. The second column tabulates the
frequencies of the peaks and the third column tabulates the harmonics
$z_1A+Z_2B+Z_3C$ of the frequencies $A$, $B$ and $C$ indicated on the
frequency axis of the power spectrum.}
\end{figure}

Figure~\rf{1000} shows the results of a simulation of the coupled rod
system for the perturbation of~\rf{900} corresponding to the
values~\rf{1001}. In the simulation the normalized rotation vector
$\bar\tau$ moved on the unit sphere as shown in the top left of the
Figure; the sense of the rotation is clockwise as seen from the top,
so that $\eta_e$ is pointing down and away from you as you look upon
the Figure. The height of the rotation vector has been graphed in the
top right of Figure~\rf{1000}, and is visibly quasiperiodic with about
two frequencies. The measured average height of the motion was
$.4970$, compared with the prediction of $.5300$. In the
middle of Figure~\rf{1000} is the power spectrum of one of the
horizontal components of $\tau$. As is shown in the table immediately
below the power spectrum, the power spectrum is consistent with the 3
fundamental angular frequencies
\begin{equation}\lb{2005} .04449,\quad .05791,\quad .1687,
\end{equation}
all other peaks being harmonics of these. The first two of these are
by far the largest peaks (note the vertical logarithmic scale in the
power spectrum) and they agree with the the predicted
frequencies~\rf{2004}. The third observed frequency in~\rf{2005}
corresponds to a far less prominent peak and is likely a frequency
associated to the higher order terms that have been truncated in the
drift approximation.

To give an idea of the speed of the drift, on the short simulation
corresponding to the top left and right of Figure~\rf{1000}, the
system rotated about 300 times while $\bar\tau$ rotated about 
7~times. Thus the drift was about 50 times slower than the original
relative equilibrium. Due to this, these kinds of simulations can be
long time: over the duration of the simulation that generated the
power spectrum in Figure~\rf{2} the rods rotated about 10,000 times.

\subsubsection{Third comparison: a singularity-induced phase jump.}

As shown in Section~(4.2), the singularities of the $-j_1=j_2$ reduced
phase spaces of the drift system are occupied by an unstable
equilibria corresponding to $SO(3)$~relative equilibria of the drift
system. A single homoclinic orbit, tracing out the intersection of
each reduced space with the plane $w_2=0$, emanates from every such
equilibrium. These homoclinic orbits are interesting features of the
drift system, in part because they provide an opportunity to
investigate dynamics near singularities of reduced spaces. I will show
here that, due to the singularity of the reduced space, the
drift system suffers a jump in its reconstruction phase as initial
conditions traverse the homoclinic orbits, and (numerically) that this
feature of the drift system persists to the coupled rod system.
The reconstruction phase jump
shows  the presence of Hamiltonian monodromy in the completely
integrable drift system; for more information on Hamiltonian
monodromy see~\ct{CushmanRHBatesLM-1997.1}, page~175, the summary
on page~403 of the same reference, as well as~\ct{BatesLM-1991.1}.

So choose $j_1$ small and consider perturbations of the relative
equilibrium
\begin{equation*}
\pi=j_1\mathbf k,\quad x=0,\quad A=\Id,
\end{equation*}
where by abuse of notation $j_1$ serves as the constant value of the
Casimir~\rf{70} of the same name. This relative equilibrium corresponds to
the equilibrium
\begin{equation}
\lb{550}w_1=w_2=0,\quad \pi_3=j_1
\end{equation}
on the $SO(3)\times S^1$~reduced space. Near~\rf{550} the motion on
the reduced  space of the  drift system  is periodic along the
intersection of the $w_2=h$ plane where $h$ is small: for a long time
the system remains near the singularity, and then moves off, passes
near $\pi_3=j_1$ and then returns to the singularity, similar to, for
example, the motion of an inverted pendulum. One $S^1$ reconstruction
gives the motion  on the phase space $\{(x,\pi)\}$ and, since the
motion is periodic on the $SO(3)\times S^1$~reduced space, there is a
well defined \defemph{$S^1$ reconstruction phase}. On the phase space
$\{(\pi,x)\}$, the only points on the level sets of $j_1$ and $j_2$
that map to~\rf{550} are
\begin{equation}\lb{551}
\pi=j_1\mathbf k,\quad x=0,
\end{equation}
since $\pi_3=j_1=-j_2$ implies by~\rf{70} and~\rf{71} that
$\pi_1=\pi_2=x_1=x_2=0$. Thus, when the $SO(3)\times S^1$ reduced
system is near~\rf{550}, the reduced system with phase
space~$\{(x,\pi)\}$ is near~\rf{551}. An $SO(3)$
reconstruction gives the motion on the phase space $\mathbb R\times
T^*SO(3)=\{(x,A,\pi)\}$, starting, say, at $A=\Id$. By
conservation of the momentum $J=A\pi$, the variable $A$ is nearly a
rotation about $\mathbf k$ whenever the reduced system is
near~\rf{551}. Thus, the motion of the point
$A\mathbf k$ on the unit sphere is this: \emph{for a long time
$A\mathbf k$ remains near $\mathbf k$, then moves off, passes near
$-\mathbf k$, and then returns near to $\mathbf k$. This motion
repeats but rotated (with respect to the previous excursion from
$\mathbf k$) by a excursion-independent angle about $\mathbf k$. This
rotation is what I mean by the \defemph{$SO(3)$~reconstruction
phase shift} of the motion.}

\begin{figure}\setlength{\unitlength}{1in}\centering
\begin{picture}(4.5,2)
\put(1.,1){\makebox(0,0){\epsfbox{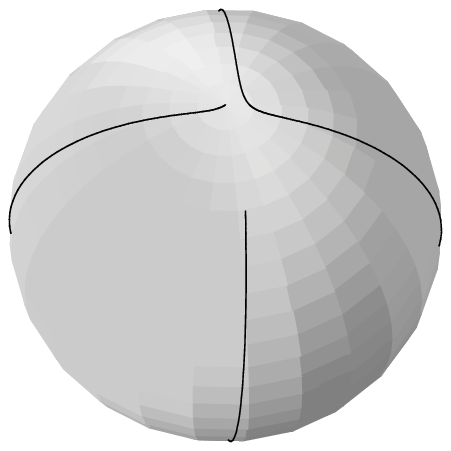}}}
\put(3.4,1){\makebox(0,0){\epsfbox{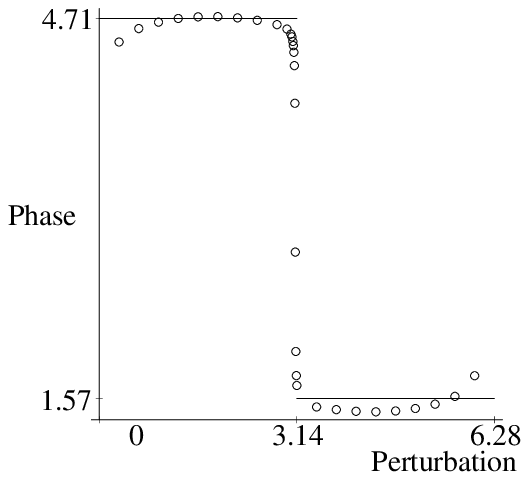}}}
\end{picture}
\caption{\lb{4002}\protect\footnotesize\it Left: the motion of $\bar\tau$
for a perturbation corresponding to being near the homoclinic orbit of
the reduced space of the drift system. The vector $\bar\tau$ begins
near the top of the sphere and moves off to the left, comes around the
sphere, makes an approximate $\pi/2$ turn, then moves away from you,
and goes around the sphere again. On the right: the $SO(3)$
reconstruction phase of the coupled rod system as initial conditions
are rotated around the singularity of the reduced phase space of the
drift system.}
\end{figure}
By simulation of the coupled rod system I have verified that the gross
details of the motion of $A\mathbf k$ in the drift system also occur
in the coupled rod system (see the left of Figure~\rf{4002}). This
situation is robust; no particular care is required, in choosing the
coupled rod system's initial conditions to evoke these kinds of
motions. The robust nature under addition of higher order terms of the
homoclinic orbit is expected through averaging theorems such as the
one in~\ct{GuckenheimerJHolmesP-1983.1},~page~168.

It is the behavior of the reconstruction phase of the drift system
near to $\pi=0$ that is relevant for predictions of the motion of the
coupled rod system near to its resonant relative equilibrium. One way
to examine this behavior is to numerically integrate initial
conditions of the drift system that are reconstructions of initial
conditions starting at a small circle surrounding the singularity of
the $SO(3)\times S^1$ reduced space. For the $S^1$ reconstruction
phase shift the result is a constant phase shift of $\pm\pi/2$ with a
jump of $\pi$ as the homoclinic orbit is traversed. The corresponding
phase shift of $A\mathbf k$ in fact becomes undefined as the
perturbation vanishes. The reason for this is that 
the equation of motion for $A$, namely Equation~\rf{55},
becomes, near to~\rf{551}, the equation
\begin{equation*}
\frac{dA}{dt}=(I_2j_1+a)\mathbf k^\wedge,
\end{equation*}
while $A(t)\mathbf k$ is not exactly $\mathbf k$ due the the presence
of the perturbation itself. Consequently the $A\mathbf k$ picks up a
rotation of angular frequency $I_2j_1+a$ during its long visit
of~\rf{551}. As the perturbation vanishes this long visit becomes an
eternity and this undefines the $SO(3)$ reconstruction phase shift as an
asymptotic effect. However, viewed from a frame that corotates with
the same angular frequency, namely $I_2j_1+a$, one can expect a well
defined phase shift. I have numerically verified this: from the
corotating frame the $SO(3)$ reconstruction phase shift of the drift
system is nearly $\pm\pi/2$ with a jump of $\pi$ as the homoclinic
orbit is traversed. This compares favorably with the $SO(3)$
reconstruction phase shift of the coupled rod system shown in right of
Figure~\rf{4002}.

Finally, so that the comparison of these phases of the drift system
and the coupled rod system are not entirely numerical, I give a
calculation showing that, asymptotically as the homoclinic orbit is
approached, these phases arise mostly from the singularity of the
reduced space. The $S^1$ reconstruction phase of the drift system may
be calculated by a slight modification (to allow reparameterization)
of the usual reconstruction method found
in~\ct{AbrahamRMarsdenJE-1978.1}: Generally, suppose $H$ is a
Hamiltonian on a symplectic phase space $P$ and that a curve
$c_\mu(s)$ is a reparameterization an evolution $c_\mu(t)$ on some
Marsden-Weinstien reduced phase
space~$P_\mu=J^{-1}(\mu)/G_\mu$. Choose a curve $d(s)\in J^{-1}(\mu)$
such that $c=\pi_\mu\circ d$, where $\pi_\mu:J^{-1}(\mu)\rightarrow
P_\mu$ is the quotient projection. Then there is a unique smooth curve
$\xi(s)\in\mathfrak g_\mu$ and a unique smooth function $a(s)$ such
that
\begin{equation}\lb{554}
X_H\bigl(d(s)\bigr)=\xi(s)p(s)+a(s)d^\prime(s).
\end{equation}
By differentiation, the curve $g\bigl(s(t)\bigr)d\bigl(s(t)\bigr)$
satisfies Hamilton's equations on $P$ if
\begin{equation*}
g^{-1}\frac{dg}{ds}=\frac1{a(s)}\xi(s),\quad\frac{dt}{ds}=\frac1{a(s)}.
\end{equation*}
Perfectly obvious generalizations hold if $d$ is just a $1$-manifold
in $J^{-1}(\mu)$ covering the image of $c_\mu$. Particularly, by
setting $\pi_2=0$, $d$ can be the subset of $\Pdrift$ defined by
\begin{gather}
\pi_1^2+\pi_2^2=j_1^2,\quad \pi_2=0,\quad
x_1=\frac{h}{2\pi_1},\quad x_2=\frac{-w_1}{2\pi_1},\lb{5553}\\
w_1^2+h^2=8(\pi_3^2-j_1^2)(\pi_3-j_1),\lb{553}
\end{gather}
in which case one calculates using~\rf{54},~\rf{554},~\rf{5553} 
and~\rf{553} that
\begin{equation}\lb{555}
\frac{d\theta}{d\pi_3}=\frac{-h\pi_3}{({j_1}^2-{\pi_3}^2)w_1},\quad
\frac{d\pi_3}{dt}=-\frac{\kappa w_1}2
\end{equation}
where $\theta$ is the $S^1$ reconstruction phase as the
counterclockwise angle of the vector $\bigl(\pi_1(t),\pi_2(t)\bigr)$
from the $\pi_2=0$ axis.

Let $r_1(h)<r_2(h)<r_3(h)$ be the roots of the cubic
$8(\pi_3-j_1)^2(\pi_3+j_1)-h^2$. By standard perturbation arguments
\begin{gather*}
r_1=-j_1+\frac1{32{j_1}^2}h^2+O(h^3),\\
r_2=j_1-\frac1{4\sqrt{j_1}}h+O(h^2),\\
r_3=j_1+\frac1{4\sqrt{j_1}}h+O(h^2).
\end{gather*}
The $S^1$ reconstruction phase over the curve
$w_1^2+h^2=8(\pi_3^2-j_1^2)(\pi_3-j_1)$ is calculated as follows. Let
the reduced system start at $w_1=0$, $\pi_3=r_2$ at time $t=t_0$, move to
$w_1=0$, $\pi_3=r_1$ at time $t=t_1$ and then complete its periodic
orbit by moving back to $w_1=0$, $\pi_3=r_2$ at time $t_2$. Then the
contribution of~\rf{555} to phase shift over the interval $[t_0,t_1]$
is
\begin{equation*}
\phi(t_1)-\phi(t_0)=\int_{t_0}^{t_1}\frac{d\phi}{dt}\,dt
=\int_{\pi_3(t_0)}^{\pi_3(t_1)}\frac{d\phi}{dt}\frac{dt}{d\pi_3}\,d\pi_3
=\int_{r_2}^{r_1}\frac{d\phi}{d\pi_3}\,d\pi_3,
\end{equation*}
where the positive square root must be used when solving for $w_1$
in~\rf{553}, since $\pi_3$ must immediately decrease after time
$t=t_0$, and by the second of~\rf{555}, $w_1$ is positive over the
interval $[t_0,t_1]$. The contribution of~\rf{555} over the interval
$[t_1,t_2]$ is identical (use the negative root here). Thus, the
total phase shift $\phi\equiv\phi(t_2)-\phi(t_0)$ over the loop as the
loop approaches the homoclinic orbit through $h>0$ is
\begin{equation}\lb{557}
\phi=-2\lim_{h\rightarrow0^+}\int_{r_1}^{r_2}\frac{-h\pi_3}{({j_1}^2-{\pi_3}^2)
\sqrt{8(\pi_3-r_1)(r_2-\pi_3)(r_3-\pi_3)}}\,d\pi_3.
\end{equation}
For small $h$ the integrand of~\rf{557} is small away from its two
singularities at $\pi_3=r_1$ and $\pi_3=r_2$. I begin with the left
singularity at $\pi_3=r_1$, so I calculate
\begin{equation*}
\phi_1\equiv\lim_{h\rightarrow0^+}\int_{r_1}^0
 \frac{-h\pi_3}{({j_1}^2-{\pi_3}^2)
 \sqrt{(\pi_3-r_1)(r_2-\pi_3)(r_3-\pi_3)}}\,d\pi_3.
\end{equation*}
Elementary estimates show that zero error is made
as $h\rightarrow0$ by the replacement of
\begin{equation*}
\frac{\pi_3}{(j_1-\pi_3)\sqrt{(\pi_3-r_2)(\pi_3-r_3)}}
\end{equation*}
with its limit as $h\rightarrow0$ of its evaluation at $\pi_3=r_1$, which is
\begin{equation*}
\lim_{h\rightarrow0^+}\frac{r_1}{(j_1-r_1)\sqrt{(r_2-r_1)(r_3-r_1)}}=
\frac{-j_1}{2j_1(\sqrt{2j_1})^2}=\frac{-1}{4j_1}.
\end{equation*}
Thus
\begin{equation*}\begin{split}
\phi_1=&\lim_{h\rightarrow0^+}
 \frac {-h}{4j_1\sqrt 8}\int_{r_1}^0
 \frac{-1}{(j_1+\pi_3)\sqrt{\pi_3-r_1}}\,d\pi_3\\
=&\lim_{h\rightarrow0^+}\frac h{4j_1\sqrt 8}
 \frac 2{\sqrt{j_1+r_1}}\arctan
 \left.\left(\frac{\sqrt{\pi_3-r_1}}{\sqrt{j_1+r_1}}\right)\right|_{\pi_3=r_1}
 ^{\pi_3=0}\\
=&\frac h{4j_1\sqrt 8}
 \frac2{\sqrt{\frac{h^2}{32{j_1}^2}}}\frac\pi2\\
=&\frac\pi2.
\end{split}\end{equation*}
The right singularity of~\rf{557} arises corresponds to the singularity of
the reduced spaces of the drift system, and it gives a phase shift of
\begin{equation*}\begin{split}
\phi_2&\equiv
 \lim_{h\rightarrow0^+}\int_0^{r_2}\frac{-h\pi_3}{({j_1}^2-{\pi_3}^2)
 \sqrt{8(\pi_3-r_1)(r_2-\pi_3)(r_3-\pi_3)}}\,d\pi_3\\
&=\frac{-h}{\sqrt8}
 \lim_{h\rightarrow0^+}\frac{\pi_3}{(j_1+r_2)\sqrt{r_2-r_1}}
 \lim_{h\rightarrow0^+}\int_0^{r_2}\frac{1}{(j_1-\pi_3)
 \sqrt{(r_2-\pi_3)(r_3-\pi_3)}}\,d\pi_3\\
&=\frac{-h}{\sqrt8}\frac{j_1}{2j_1\sqrt{2j_1}}
 \lim_{h\rightarrow0^+}\Biggl(\frac{-1}{\sqrt{(r_3-j_1)(j_1-r_2)}}\\
&\qquad\times\left.\arctan\left(
 \frac12\frac{\pi_3(2j_1-r_2-r_3)+2r_2r_3-j_1(r_2+r_3)}
 {\sqrt{(r_3-j_1)(j_1-r_2)(r_2-\pi_3)(r_3-\pi_3)}}\right)
 \right|_{\pi_3=0}^{\pi_3=r_2}\Biggr)\\
&=\frac{-h}{\sqrt8}\frac{j_1}{2j_1\sqrt{2j_1}}
 \frac{-1}{\frac{h}{4\sqrt{j_1}}}\frac{-\pi}{2}\\
&=-\frac\pi4,
\end{split}\end{equation*}
so that the phase shift as the homoclinic orbit is approached through
positive $h$ is
\begin{equation*}
\phi=-2(\phi_1+\phi_2)=-2\left(\frac\pi2-\frac\pi4\right)=-\frac\pi2.
\end{equation*}
The same calculation for $h$ negative gives $\phi_1=-\pi/2$ and
$\phi_2=\pi/4$ for a total phase shift of $\phi=\pi/2$.

\footnotesize\frenchspacing

\end{document}

%% file: fig1.tex
\begin{picture}(0,0)%
\epsfbox{fig1.pstex}%
\end{picture}%
\setlength{\unitlength}{0.00083300in}%
\begin{picture}(5456,6420)(69,-5706)
\put(3571,-3176){\makebox(0,0)[lb]{\scriptsize$\pi_3$}}
\put(4937,-1307){\makebox(0,0)[lb]{\scriptsize$\pi_3$}}
\put(2931,-5351){\makebox(0,0)[lb]{\scriptsize$-j_1$}}
\put(2036,-2311){\makebox(0,0)[lb]{\scriptsize$j_1$}}
\put(361,-811){\makebox(0,0)[lb]{\scriptsize$j_2$}}
\put(3168,-572){\makebox(0,0)[lb]{\scriptsize$j_1$}}
\put(2975,-3204){\makebox(0,0)[lb]{\scriptsize$-j_2=j_1$}}
\put(3309,-3634){\makebox(0,0)[lb]{\scriptsize$-j_1$}}
\put(4657,-5501){\makebox(0,0)[lb]{\scriptsize$\pi_3$}}
\put(121,-4382){\makebox(0,0)[lb]{\scriptsize$j_1=0$}}
\put(1036,123){\makebox(0,0)[lb]{\scriptsize$j_1$}}
\put(1287,-25){\makebox(0,0)[lb]{\scriptsize$\pi_3$}}
\put(2594,-761){\makebox(0,0)[lb]{\scriptsize$-j_2$}}
\put(2192,-586){\makebox(0,0)[lb]{\scriptsize$-j_1$}}
\put(3602,-736){\makebox(0,0)[lb]{\scriptsize$\pi_3$}}
\put(3256,-1346){\makebox(0,0)[lb]{\scriptsize$-j_1$}}
\put(4543,-1163){\makebox(0,0)[lb]{\scriptsize$-j_2$}}
\put(2347,-3021){\makebox(0,0)[lb]{\scriptsize$-j_1$}}
\put(5019,-3783){\makebox(0,0)[lb]{\scriptsize$\pi_3$}}
\put(4598,-3631){\makebox(0,0)[lb]{\scriptsize$-j_2=j_1$}}
\put(3248,-4770){\makebox(0,0)[lb]{\scriptsize$\pi_3$}}
\put(4262,-5526){\makebox(0,0)[lb]{\scriptsize$j_1$}}
\put(1014,-4382){\makebox(0,0)[lb]{\scriptsize$-j_2$}}
\put(1324,-4530){\makebox(0,0)[lb]{\scriptsize$\pi_3$}}
\put(210,123){\makebox(0,0)[lb]{\scriptsize$j_2=-j_1$}}
\put(2574,-4808){\makebox(0,0)[lb]{\scriptsize$j_1$}}
\put(2890,-4635){\makebox(0,0)[lb]{\scriptsize$-j_2$}}
\put(2031,-4635){\makebox(0,0)[lb]{\scriptsize$-j_1$}}
\end{picture}